\documentclass[final]{siamart171218}

\usepackage{amsfonts,bm}
\usepackage{graphicx}
\usepackage{epstopdf}
\usepackage{algorithmic}
\usepackage{cases}
\usepackage{subfigure,multirow}
\usepackage{amsopn}

\makeatletter
\@addtoreset{Equation}{Section}
\makeatother

\theoremstyle{remark}
\newtheorem{remark}[theorem]{Remark}

\newcommand{\TheTitle}{A numerical approach to the optimal control of thermally convective flows}
\newcommand{\TheAuthors}{Yongcun Song, Xiaoming Yuan, Hangrui Yue}

\headers{A numerical approach to the optimal control of thermally convective flows}{Y. Song, X. Yuan, H. Yue}

\title{{\TheTitle}\thanks{28 November, 2022
		\funding{The work of the first author was supported by the Humboldt Research Fellowship for postdoctoral researchers. The work of the second author was supported by a seed fund for basic research at The University of Hong Kong. The work of the third author was supported by the Fundamental Research Funds for the Central Universities, Nankai University (Grant Number 63221035) }}}

\author{
Yongcun Song\thanks{Chair for Dynamics, Control and Numerics, Alexander von Humboldt-Professorship, Department of Data Science,  Friedrich-Alexander-Universit\"at Erlangen-N\"urnberg, 91058 Erlangen, Germany.  \email{ysong307$@$gmail.com.}}
\and
Xiaoming Yuan\thanks{Department of Mathematics, The University of Hong Kong. Hong Kong, China. \email{xmyuan@hku.hk}}
\and
Hangrui Yue\thanks{School of Mathematical Sciences, Nankai University, Tianjin, China. \email{yuehangrui@gmail.com}}
}

\ifpdf
\hypersetup{
	pdftitle={\TheTitle},
	pdfauthor={\TheAuthors}
}
\fi

\begin{document}
	
\maketitle
\begin{center}
	{\em  Dedicated to the memory of Roland Glowinski (1937--2022): a dear friend and great mentor.\footnote{The authors were encouraged by Roland Glowinski to consider this work when he was visiting Hong Kong in 2019. }}

\end{center}

\bigskip
\begin{abstract}
The optimal control of thermally convective flows is usually modeled by an optimization problem with constraints of Boussinesq equations that consist of the Navier-Stokes equation and an advection-diffusion equation. This optimal control problem is challenging from both theoretical analysis and algorithmic design perspectives. For example, the nonlinearity and coupling of fluid flows and energy transports prevent direct applications of gradient type algorithms in practice. In this paper, we propose an efficient numerical method to solve this problem based on the operator splitting and optimization techniques. In particular, we employ the Marchuk-Yanenko method leveraged by the $L^2-$projection for the time discretization of the Boussinesq equations so that the Boussinesq equations are decomposed into some easier linear equations without any difficulty in deriving the corresponding adjoint system. Consequently, at each iteration, four easy linear advection-diffusion equations and two degenerated Stokes equations at each time step are needed to be solved for computing a gradient. Then, we apply the Bercovier-Pironneau finite element method for space discretization, and design a BFGS type algorithm for solving the fully discretized optimal control problem. We look into the structure of the problem, and design a meticulous strategy to seek step sizes for the BFGS efficiently. Efficiency of the numerical approach is promisingly validated by the results of some preliminary numerical experiments.

\end{abstract}


%% REQUIRED

\section{Introduction}
Thermally convective flows arise in many applications such as thermal insulation, cooling of fluids in channels surrounding nuclear reactor core, the circulation of liquid metals in solidifying ingot, the manufacture of crystals, and modeling, designing, and controlling energy-efficient building systems, see e.g., \cite{B08,glowinski2003finite, SSR02}. Such flows are typically modeled by the Boussinesq equations, which consist of the Navier-Stokes equation for incompressible viscous flow coupled with an advection-diffusion equation for temperature.  In practice,
thermally convective flows also play a crucial role in the control of different physical processes. In view of this, various controls of thermally convective flows have been widely studied from different perspectives in the past decades. For instance, the optimal control of temperature peaks along the bounding surfaces of containers of fluid flows has been studied in \cite{GL94}. In \cite{BKR98}, feedback control for thermal fluid was considered, and in \cite{IR98} the optimal control of flow separation in a channel flow using temperature control was investigated. Linear quadratic regulator control and pointwise control of the Boussinesq equations were studied in \cite{BH13} and \cite{NR11}, respectively.  More related literature on the control of thermally convective flows can be referred to \cite{abergel1990some,B08, BH06,belmiloudi2000,belmiloudi2002, IR98,casas1994, Lee2003,LI00, LI2000, SSR97, SSR02} and references therein.

\subsection{Problem statement and motivation}

 The class of thermally convective flow under our consideration is modeled by the following non-dimensional Boussinesq equations: \begin{flalign}\label{state_equation1}
 	\left\{
 	\begin{aligned}
 		&\frac{\partial y}{\partial t}+(y\cdot{\nabla})y-\nu_1 \Delta y+\nabla p=\theta \mathbf{{e}}_2\quad \text{in}\quad Q, \\
 		&\nabla\cdot y=0 \quad \text{in}\quad Q,\\
 		&\frac{\partial \theta}{\partial t}-\nu_2 \Delta\theta+y\cdot\nabla\theta=0\quad \text{in}\quad Q,
 	\end{aligned}
 	\right.
 \end{flalign}
 together with the boundary and initial conditions
 \begin{flalign}\label{state_equation2}
 \left\{
 	\begin{aligned}
 		&y=0&~ \text{on} ~\Sigma,\quad
 		y(0)=y_0.&\\
 		&\nu_2\frac{\partial\theta}{\partial \vec{n}}=v\chi_{\Sigma_c}& ~\text{on} ~\Sigma,\quad
 		\theta(0)=\theta_0.&
 	\end{aligned}
 \right.
 \end{flalign}
Above, $Q=\Omega\times (0,T)$ and  $\Sigma=\Gamma\times (0,T)$, where $\Omega\subset \mathbb{R}^2$ is a bounded domain with Lipschitz continuous boundary $\Gamma=\partial\Omega$.  The function $\chi_{\Sigma_c}$ is the characteristic function of $\Sigma_c=\Gamma_c\times (0,T)$ with $\Gamma_c\subset\Gamma$. The variable $y=(y_1,y_2)^T$ is the flow velocity, $p$ is the pressure deviation from the hydrostatic and $\theta$ is a normalized temperature deviation. The coefficients $\nu_1=\sqrt{\frac{Pr}{Ra}}$ and $\nu_2=\frac{1}{\sqrt{RaPr}}$ are constants with $Pr$ and $Ra$ the Prandtl and Rayleigh numbers, respectively. The vector $\mathbf{{e}}_2=(0,1)^T$, and $\vec{n}$ is the unit outward normal vector. Equation (\ref{state_equation1})-(\ref{state_equation2})  describes thermal convection, with the flow induced by gravity and differences in temperature on the boundary $\Gamma$. For brevity, we focus on the boundary conditions \eqref{state_equation2} in the following discussions and other types of boundary conditions can be treated similarly.

 We study the optimal control of thermally convective flows modeled by (\ref{state_equation1})-(\ref{state_equation2}) via the heat flux $v\in L^2(\Sigma_c)$, which can be mathematically expressed as
\begin{flalign}\label{OptimalControl}
	\left\{
	\begin{array}{ll}
		u\in L^2(\Sigma_c), \\
		J(u)\leq J(v), \forall v\in L^2(\Sigma_c),
	\end{array}
	\right.
\end{flalign}
where the objective functional $J$ is given by
\begin{equation}\label{objective_functional}
	J(v)=\frac{1}{2}\iint_Q|y-y_d|^2dxdt+\frac{\alpha}{2}\iint_{\Sigma_c} v^2dxdt,  \forall v\in L^2(\Sigma_c),
\end{equation}
 the constant $\alpha>0$ is a regularization parameter, and the state $y=y(t;v)$ obtained from the control $v$ is the solution to the equation (\ref{state_equation1})-(\ref{state_equation2}).
The optimal control problem (\ref{OptimalControl})--(\ref{objective_functional}) aims at determining the velocity-pressure-temperature triplet $(y, p, \theta)$-- by controlling the heat flux $v$ -- such that $y$ best matches a target velocity field $y_d$. Such a model plays an important role in many real applications, e.g., optimal control of semiconductor melts in zone-melting and Czochralski growth configurations \cite{BH06}, and control of energy efficient building systems  \cite{bbsz09}.  Therefore, it is of practical significance to develop an efficient numerical method for solving  (\ref{OptimalControl})--(\ref{objective_functional}).

 %{\color{red} Actually, investigating this kind of problem numerically was suggested by R. Glowinski when he visited Hong Kong in 2019. }

In addition to the tracking-type objective functional given by (\ref{objective_functional}), one may also be interested in the following vorticity reduction objective functional
\begin{equation}\label{objective_functional_r}
J(v)=\frac{1}{2}\iint_Q|\nabla\times y(v)|^2dxdt+\frac{\alpha}{2}\iint_{\Sigma_c} v^2dxdt,  \forall v\in L^2(\Sigma_c),
\end{equation}
where $\frac{1}{2}\iint_Q|\nabla\times y(v)|^2dxdt$ measures the vorticity of the flow. Minimizing (\ref{objective_functional_r})  has significant applications in science and engineering
such as control of turbulence and control of crystal growth process \cite{B08, BH06}. To expose our main ideas clearly, we focus on the objective functional (\ref{objective_functional}) and all results can be easily applied to the case of (\ref{objective_functional_r}) by slight modifications.

In the equation (\ref{state_equation1})-(\ref{state_equation2}), the Navier-Stokes equation and an advection-diffusion are coupled; the nonlinearity of $(y\cdot{\nabla})y$ and $(y\cdot{\nabla})\theta$ are coupled with the incompressibility condition $\nabla\cdot y=0$; and the problem (\ref{OptimalControl})--(\ref{objective_functional}) is highly nonconvex because of the nonlinearity.
The complicated structure and nonconvexity lead to enormous obstacles in solving (\ref{OptimalControl})--(\ref{objective_functional}). In particular, as to be shown in Theorem \ref{thm:oc}, computing a gradient of  $J(v)$ requires solving equation (\ref{state_equation1})-(\ref{state_equation2}) and the corresponding adjoint system, and thus is computationally challenging. Consequently, although gradient-type methods (e.g., gradient descent methods and conjugate gradient methods) can be conceptually applied to  (\ref{OptimalControl})--(\ref{objective_functional}),  it is rather nontrivial to implement them in practice. All these obvious difficulties explain that research on numerical study for (\ref{OptimalControl})--(\ref{objective_functional}) is still in its infancy.

Some easier cases of (\ref{OptimalControl})--(\ref{objective_functional})  have been studied in \cite{Lee2003,LI00,LI2000}, where the state equation (\ref{state_equation1})-(\ref{state_equation2}) was replaced by its stationary counterpart. In \cite{HRH2021}, robust temperature and velocity output tracking for linearized Boussinesq equations was studied. For the general nonlinear time-dependent case (\ref{state_equation1})-(\ref{state_equation2}),  the exact controllability was analyzed in \cite{FI99}, the existence of a solution and the first-order optimality condition was studied in \cite{BH06,barwolff2007analysis,casas1994}. In \cite{BHH2016}, local exponential stabilization of (\ref{state_equation1})-(\ref{state_equation2}) was studied with control acting on {\color{blue} a} portion of the boundary. Computationally,  a semi-implicit scheme was suggested in \cite{B08, BH06} for the time discretization and a damped Picard iteration was proposed for solving the first-order optimality system of (\ref{OptimalControl})--(\ref{objective_functional}). Note that the semi-implicit time discretization scheme is only conditionally stable, and in practice, a tiny time step might be necessary to ensure numerical stability and robustness. The damped Picard iteration is essentially a gradient descent method, which usually converges slowly.  In \cite{SSR97}, a piecewise-in-time optimal control approach was proposed,  which matches the velocity fields $y_d$ at a sequence of time intervals, and the velocity tracking at each time interval is then formulated as an optimal control problem modeled by a stationary Boussinesq equation. A solution of (\ref{OptimalControl})--(\ref{objective_functional}) was then obtained by patching together all the solutions of local optimal control problems at each time interval. However, the piecewise-in-time optimal approach can only pursue a suboptimal control of (\ref{OptimalControl})--(\ref{objective_functional}). In \cite{SSR02}, an adaptive procedure was proposed by using a proper orthogonal decomposition approach with the sequential quadratic programming method to obtain a reduced-order model for the optimal control problem (\ref{OptimalControl}) with the vorticity reduction objective functional (\ref{objective_functional_r}).
Note that the backward Euler method was used in \cite{SSR02} for the time discretization, which leads to a large-scale and computationally expensive complex nonlinear system at each time step with all components coupled together.

\subsection{Our methodology}

In this paper, we propose an efficient and relatively easy-to-implement numerical method to solve (\ref{OptimalControl})--(\ref{objective_functional}). We employ the discretize-then-optimize approach, i.e., we first discretize the optimal control problem (\ref{OptimalControl})--(\ref{objective_functional}) and then compute the gradient in a fully discrete setting. Compared with the optimize-then-discretize approach, the discretize-then-optimize approach is more advantageous in that the fully discrete state equation and the fully discrete adjoint equation are strict in duality. This fact guarantees that the computed negative gradient is a descent direction of the fully discrete optimal control problem. We refer to, e.g. \cite{GH1998,zuazua2015}, for more discussions on the difference between the discretize-then-optimize and optimize-then-discretize approaches. To implement the discretize-then-optimize approach, we advocate combining an operator splitting scheme for time discretization and a finite element method for space discretization. Then, the Broyden-Fletcher-Goldfarb-Shanno (BFGS) quasi-Newton method is applied to solve the resulting optimization problem. Although these numerical approaches have been individually developed in different literatures, it is highly nontrivial to implement them synergically to solve the problem (\ref{OptimalControl})--(\ref{objective_functional}) due to its complicated structure.

Recall that when a gradient-type method is applied to solve (\ref{OptimalControl})--(\ref{objective_functional}), computing the gradient of the objective functional $J(v)$ requires solving the state equation (\ref{state_equation1})--(\ref{state_equation2}) and the corresponding adjoint equation, which is usually challenging and computationally expensive in numerical implementation. To address this issue, we aim at developing an easily implementable time discretization scheme for (\ref{OptimalControl})--(\ref{objective_functional}) such that the resulting discrete gradient is cheap to compute. Our philosophy is to decompose the state equation (\ref{state_equation1})--(\ref{state_equation2}) into smaller and easier parts so that the computation of a discrete gradient only requires solving several much easier subproblems. More precisely, we aim at untieing the velocity-pressure-temperature triplet, and decoupling the nonlinearity of $(y\cdot{\nabla})y$ and $(y\cdot{\nabla})\theta$ from the incompressibility condition $\nabla\cdot y=0$. For this purpose, we consider the operator splitting techniques. In the literature, it has been shown that various operator splitting methods such as the Douglas-Rachford splitting method in \cite{douglas1956}, the Peaceman-Rachford splitting method in \cite{peaceman1955}, and the $\theta$-scheme in \cite{glowinski1985,glowinski1986} are very efficient for solving various partial differential equations (PDEs)  (see \cite{glowinski2003finite}). Their applications to optimal control problems, however, have not been well explored. It is worth noting that a straightforward application of the aforementioned operator splitting methods to the optimal control problem (\ref{OptimalControl})--(\ref{objective_functional}) leads to some immediate difficulties in numerical implementation. In particular, as discussed in Remark \ref{br1},  deriving the corresponding adjoint equation is challenging due to some coupled terms at different time intervals. This issue has also been mentioned in \cite{berggren1998}, where a semi-implicit finite difference scheme was used for the time discretization of an optimal control problem modeled by the Navier-Stokes equation. The author commented ``There are other more sophisticated schemes that have better accuracy and stability properties. Examples are the operator-splitting schemes" and ``However, such schemes will introduce some nontrivial complications when deriving corresponding adjoint equation". To the best of our knowledge, there has been almost no development on implementing operator splitting methods to optimal control problems including (\ref{OptimalControl})--(\ref{objective_functional}).

To tackle the aforementioned issues, we propose to leverage the Marchuk-Yanenko method \cite{Marchuk1990} and the $L^2$--projection \cite{Chorin} to implement the time discretization of equation (\ref{state_equation1})-(\ref{state_equation2}). The Marchuk-Yanenko method is chosen because it does not lead to coupled terms at different time intervals and hence does not introduce any complication in deriving the corresponding adjoint equation. We are motivated to consider the $L^2$-projection method by its popular application for unsteady incompressible Navier-Stokes equations, see \cite[Chapter 7]{glowinski2003finite} and the survey paper \cite{GMS05} for a thorough discussion. Inspired by \cite{GMS05},  we advocate the $L^2$-projection with an incremental term of the pressure $p$ (see (\ref{deg_Stokes_y})) to increase numerical accuracy and stability. Consequently, a scheme combining the Marchuk-Yanenko method with the incremental $L^2$-projection is proposed for the time discretization of (\ref{state_equation1})-(\ref{state_equation2}).  The resulting scheme only needs to solve a sequence of decoupled linear time-independent equations for $y$, $p$, and $\theta$ at each time step. Computing time can thus be substantially lowered for large-scale cases. With the proposed time-discretization scheme, the gradient is relatively easier to compute when a gradient-type method is applied. More precisely, we only need to solve four linear advection-diffusion equations and two degenerated Stokes equations at each time step to compute the gradient.  All these equations can be easily solved by some well-developed numerical methods in the literature, e.g., the fixed-point iterative schemes in \cite{BN99} for advection-diffusion equations and the preconditioned conjugate gradient methods (e.g.,  \cite[Chapter 3]{glowinski2003finite}) for degenerated Stokes equations, respectively.

With the well-designed operator splitting time-discretization scheme, the gradient of the objective functional $J(v)$ at $v\in L^2(\Sigma_c)$ is easy to compute and thus classic gradient descent methods can be applied to solve (\ref{OptimalControl})-(\ref{objective_functional}). However, gradient descent methods may converge slowly and inefficiently. To address this issue, we suggest using the BFGS method for solving (\ref{OptimalControl})-(\ref{objective_functional}). It is well known that the BFGS method may be very effective to deal with large-scale optimization problems, see, e.g. \cite{nocedal2006sequential}. To implement the BFGS, it is crucial to seek a suitable step size for each iteration while some commonly used backtracking line search strategies are usually too expensive for the settings under our consideration, because evaluating the objective function value $J(v)$ is required repeatedly and evaluating $J(v)$ entails solving the state equation (\ref{state_equation1})--(\ref{state_equation2}). To tackle this difficulty, we look into the structure of the problem (\ref{OptimalControl})-(\ref{objective_functional}) meticulously and propose an efficient inexact strategy for determining step sizes that requires solving only a few linear equations. Thus, the implementation of the BFGS is easy and the computation for solving (\ref{OptimalControl})--(\ref{objective_functional}) becomes much cheaper.

Finally, we mention that the central concern of this work is to design an efficient numerical approach to the optimal control problem (\ref{OptimalControl})--(\ref{objective_functional}). Other techniques developed in the literatures such as the model reduction technique \cite{SSR02} and the memory-saving strategy \cite{berggren1998} can also be embedded into the implementation of our proposed numerical approach to further reduce the computational cost. Furthermore, since the Navier-Stokes equation is involved in (\ref{state_equation1})-(\ref{state_equation2}), our proposed method can be directly applied to optimal control problems modeled by the Navier-Stokes equation  \cite{berggren1998, Gunz02}.

\subsection{Organization}
The outline of this paper is as follows. In Section \ref{sec:analysis}, we present some preliminary results and derive the first-order optimality conditions for (\ref{OptimalControl})--(\ref{objective_functional}).  In Section \ref{se:time_dis}, we propose an operator splitting method for the time discretization of (\ref{OptimalControl})--(\ref{objective_functional}) and a finite element
method for the space discretization. Then, a limited-memory BFGS method with an efficient step size strategy is proposed in Section \ref{sec:Numerical algorithm} for solving the fully discretized problem of (\ref{OptimalControl})--(\ref{objective_functional}). Some preliminary numerical results are reported in Section \ref{sec:numerical_results} to validate the efficiency of
our proposed numerical approach. Finally, some conclusions are drawn in Section \ref{sec:conclusions}.

\section{Analysis of problem (\ref{OptimalControl})--(\ref{objective_functional})} \label{sec:analysis}
In this section, we introduce some notations and preliminary results, which will be used in the following discussions. Then, we show the existence of optimal controls for (\ref{OptimalControl})--(\ref{objective_functional}) and derive the associated first-order optimality conditions.
\subsection{Notations}\label{se:notations}
We shall use the standard notations for the Sobolev spaces $H^m(\Omega)$ with norm $\|\cdot\|_m$ and $L^2(\Omega)$ with norm $\|\cdot\|$. Let $H_0^m(\Omega)$ be the closure of the space $C_0^\infty(\Omega)$ under the norm $\|\cdot\|_m$, where $C_0^{\infty}(\Omega)$ denotes the space of all infinitely differentiable
functions over $\Omega$ with a compact support in $\Omega$.  Let $X$ be a Banach space with a norm $\|\cdot\|_X$. Then, the space $L^2(0, T;X)$ consists of all measurable functions $z:(0,T)\rightarrow X$ satisfying
$$
\|z\|_{L^2(0, T;X)}:=\left(\int_0^T\|z(t)\|_X^2dt \right)^{\frac{1}{2}}<+\infty.
$$
Moreover, we use the following function spaces:
\begin{gather*}
Y:=\{y\in [L^2(\Omega)]^2, \nabla \cdot y=0\},~V:=\{y\in [H_0^1(\Omega)]^2, \nabla \cdot y=0\}, ~
P:=L_0^2(\Omega)=\{p\in L^2(\Omega), \int_\Omega p~dx=0\}\\
W_y(0,T):=\{y_t\in L^2(0,T;V^*),y\in L^2(0,T;V)\},\quad
 W_\theta(0,T):=\{\theta_t\in L^2(0,T;H^1(\Omega)^*),\theta\in L^2(0,T;H^1(\Omega))\},
\end{gather*}
where $\mathcal{Z}^*$ is the dual space of $\mathcal{Z}$ with $\mathcal{Z}$ a Banach space.
To simplify the notation, we introduce the following bilinear and trilinear forms:
$$
\left\{
\begin{aligned}
&a(y,\varphi)=\nu_1\int_\Omega\nabla y\cdot\nabla \varphi dx,\quad
 \forall y,\varphi \in [H^1(\Omega)]^2,\\
 &b(q, y)=\int_\Omega q\nabla \cdot ydx, \quad q\in P, y\in [H^1(\Omega)]^2,\\
&d(\theta,\psi)=\nu_2\int_\Omega \nabla \theta\cdot\nabla \psi dx,\quad \forall \theta ,\psi \in H^1(\Omega),\\
&c(y,w,\varphi)=\int_{\Omega}(y\cdot\nabla)w\cdot \varphi~dx,\quad\forall y,w,\varphi\in[H^1(\Omega)]^2,\\
&e(y,\theta,\psi)=\int_{\Omega}(y\cdot\nabla\theta)\psi~dx, \quad\forall y\in [H^1(\Omega)]^2,\quad \theta, \psi\in H^1(\Omega).
\end{aligned}
\right.
$$
It is easy to verify that both bilinear forms $a(y,\varphi)$ and $d(\theta, \psi)$ are continuous and coercive using some standard arguments as those in e.g., \cite{IR98}.

\subsection{Existence of optimal controls}
With these notations, the variational formulation of the state equation reads as follows: find $(y, p,\theta)\in L^2(0, T;[H_0^1(\Omega)]^2)\times L^2(0, T; P)\times L^2(0, T; H^1(\Omega))$ such that for a.e. $t\in (0, T)$,
\begin{flalign}\label{variationalform}
\left\{
\begin{aligned}
&(y_t, \varphi)+a(y(t), \varphi)+c(y(t), y(t), \varphi)-b(p(t),\varphi)=(\theta(t) \bm{e}_2, \varphi)&\quad \forall \varphi\in [H_0^1(\Omega)]^2,\\
&b(q,y(t))=0&\quad \forall q\in P,\\
&(\theta_t,\psi)+d(\theta(t),\psi)+e(y(t), \theta(t),\psi)=(v(t)\chi_{\Gamma_c},\psi)_\Gamma&\quad \forall \psi\in H^1(\Omega),\\
&y(x,0)=y_0,\quad \theta(x,0)=\theta_0&\quad x\in \Omega.
\end{aligned}
\right.
\end{flalign}
Here and in what follows, the notation $(\cdot,\cdot)$ stands for the canonical $L^2$-inner product.
Moreover, if we consider the subspace $V\subset [H_0^1(\Omega)]^2$, then the variational formulation (\ref{variationalform}) can be reformulated as: find $(y,\theta)\in L^2(0,T;V)\times L^2(0,T;H^1(\Omega))$ such that for a.e. $t\in (0,T)$
\begin{flalign*}
	\left\{
	\begin{aligned}
		&(y_t, \varphi)+a(y(t), \varphi)+c(y(t), y(t), \varphi)=(\theta(t) \bm{e}_2, \varphi)&\quad \forall \varphi\in V,\\
		&(\theta_t,\psi)+d(\theta(t),\psi)+e(y(t), \theta(t),\psi)=(v(t)\chi_{\Gamma_c},\psi)_\Gamma&\quad \forall \psi\in H^1(\Omega),\\
		&y(x,0)=y_0,\quad \theta(x,0)=\theta_0&\quad x\in \Omega.
	\end{aligned}
	\right.
\end{flalign*}
For the state equation (\ref{state_equation1})$-$(\ref{state_equation2}), we have the following result on the existence and uniqueness of the solution. We refer to, such as \cite{gajewski1975iterativen}, for the details of the proof.
\begin{theorem}\label{thm:soln_state}
	Suppose that $y_0\in Y$ satisfies $y_{0}|_{\Gamma}=0$, $\theta_0\in L^2(\Omega)$ and $v\in L^2(0,T; H^{-1/2}(\Gamma))$. Then, the state equation (\ref{state_equation1})$-$(\ref{state_equation2}) admits a unique weak solution $(y,p,\theta)\in W_y(0,T)\times P\times W_\theta(0,T)$.
\end{theorem}

Using the result of Theorem \ref{thm:soln_state} and considering the continuous embeddings:
$$W_y(0,T)\hookrightarrow C(0,T;Y),$$
it can be shown that the objective functional (\ref{objective_functional}) is well-defined. The existence of a solution to the optimal control problem (\ref{OptimalControl})-(\ref{objective_functional}) can be proved following the lines in \cite{lions1971optimal}. We thus omit the proof here for succinctness, and refer to  \cite{abergel1990some,barwolff2007analysis,casas1994} for the details. Since problem (\ref{OptimalControl})-(\ref{objective_functional})  under investigation is non-convex, the uniqueness of the solution cannot be guaranteed and only a local minimizer can be pursued in general cases.

\subsection{First-order optimality conditions}
In this subsection, we use a formal perturbation analysis, which has been well developed in \cite{lions1971optimal}, to derive the first-order optimality condition of the optimal control problem (\ref{OptimalControl})-(\ref{objective_functional}).
Let $DJ(v)$ be the first-order derivative of $J$ at $v\in L^2(\Sigma_c)$. If $u\in L^2(\Sigma_c)$ is a local optimal control to (\ref{OptimalControl})-(\ref{objective_functional}),  then the optimality condition corresponding to the optimal control problem (\ref{OptimalControl})-(\ref{objective_functional}) reads as
\begin{equation*}\label{oc_sub_u}
DJ(u)=0.
\end{equation*}

To specify $DJ(v)$, let us consider $v\in L^2(\Sigma_c)$ and a small perturbation $\delta v\in L^2(\Sigma_c)$, we then have
\begin{gather}{\label{def_delta_j}}
\delta J(v)=(DJ(v), \delta v)=\iint_{Q}(y-y_d)\delta y ~dxdt+\alpha\iint_{{\Sigma_c}}v\delta v~ dxdt,
\end{gather}
where $\delta y$ is the solution of
\begin{flalign*}\label{per_state_equation1}
\left\{
\begin{aligned}
&\frac{\partial \delta y}{\partial t}+(y\cdot{\nabla})\delta y+(\delta y\cdot{\nabla})y-\nu_1 \Delta \delta y+\nabla \delta p=\delta\theta \mathbf{{e}}_2\quad \text{in}\quad \Omega\times(0,T), \\
&\nabla\cdot \delta y=0 \quad \text{in}\quad \Omega\times(0,T),\\
&\delta y=0\quad \text{on}\quad \Gamma\times(0,T),\\
&\delta y(0)=0,
\end{aligned}
\right.
\end{flalign*}
and $\delta\theta$ satisfies
\begin{gather*}\label{per_state_equation2}
\left\{
\begin{aligned}
&\frac{\partial \delta\theta}{\partial t}-\nu_2 \Delta\delta\theta+y\cdot\nabla\delta\theta+\delta y\cdot\nabla \theta=0\quad \text{in}\quad \Omega\times(0,T), \\
&\nu_2\frac{\partial\delta\theta}{\partial \vec{n}}=\delta v\chi_{\Sigma_c}\quad \text{on}\quad \Gamma\times(0,T),\\
&\delta\theta(0)=0.&
\end{aligned}
\right.
\end{gather*}

For any ${z}\in [L^2(Q)]^2$ and $\xi,\zeta\in L^2(Q)$, we have
\begin{gather}
\iint_{Q}(\frac{\partial \delta y}{\partial t}+(y\cdot{\nabla})\delta y+(\delta y\cdot{\nabla})y-\nu_1 \Delta \delta y+\nabla \delta p)\cdot {z} dxdt=\iint_{Q}\delta\theta \mathbf{{e}}_2\cdot {z} dxdt,\label{e1}\\
\iint_{Q}\nabla\cdot \delta y\xi dxdt=0, \label{e2}\\
\iint_{Q}(\frac{\partial \delta\theta}{\partial t}-\nu_2 \Delta\delta\theta+y\cdot\nabla\delta\theta+\delta y\cdot\nabla \theta)\zeta dxdt=0.\label{e3}
\end{gather}
Applying Green's formula to (\ref{e1})$-$(\ref{e3}), we obtain
\begin{gather}\label{equality1}
\begin{aligned}
\delta y(T){z}(T)-\delta y(0){z}(0)+\iint_Q(-\frac{\partial {z}}{\partial t}-\nu_1\Delta{z}+(\nabla y)^T{z}-(y\cdot \nabla){z})\cdot\delta y dxdt\\
-\iint_Q(\nabla\cdot {z})\delta p dxdt+\iint_\Sigma \frac{\partial\delta p}{\partial \vec{n}}{z} dxdt=\iint_Q\delta\theta \mathbf{{e}}_2\cdot {z} dxdt,
\end{aligned}\\
-\iint_Q \nabla \xi\cdot \delta y dxdt=0,\label{equality2}\\
\begin{aligned}
\delta \theta(T)\zeta(T)-\delta \theta(0)\zeta(0)+\iint_Q(-\frac{\partial \zeta}{\partial t}-\nu_2\Delta\zeta-y\cdot\nabla\zeta)\delta \theta dxdt\\
+\iint_Q \zeta\nabla\theta\cdot\delta y dxdt=\nu_2\iint_\Sigma (\frac{\partial \delta \theta}{\partial \vec{n}}\zeta-\frac{\partial \zeta}{\partial \vec{n}}\delta\theta)dxdt. \label{equality3}
\end{aligned}
\end{gather}
Summing up (\ref{equality1})$-$(\ref{equality3}) shows that
\begin{eqnarray}\label{sum_equality}
\begin{aligned}
&\delta y(T){z}(T)-\delta y(0){z}(0)+\delta \theta(T)\zeta(T)-\delta \theta(0)\zeta(0)\\
&+\iint_Q(-\frac{\partial {z}}{\partial t}-\nu_1\Delta{z}+(\nabla y)^T{z}-(y\cdot \nabla){z}-\nabla \xi)\cdot\delta y dxdt\\
&+\iint_Q \zeta\nabla\theta\cdot\delta y dxdt+\iint_Q(-\frac{\partial \zeta}{\partial t}-\nu_2\Delta\zeta-y\cdot\nabla\zeta)\delta \theta dxdt\\
&-\iint_Q(\nabla\cdot {z})\delta p dxdt+\iint_\Sigma \frac{\partial\delta p}{\partial \vec{n}}{z} dxdt\\
=&\iint_Q\delta\theta \mathbf{{e}}_2\cdot {z} dxdt+\nu_2\iint_\Sigma (\frac{\partial \delta \theta}{\partial \vec{n}}\zeta-\frac{\partial \zeta}{\partial \vec{n}}\delta\theta)dxdt.
\end{aligned}
\end{eqnarray}
If we let $z$ and $\zeta$ satisfy
\begin{gather*}
-\frac{\partial {z}}{\partial t}-\nu_1\Delta{z}+(\nabla y)^T{z}-(y\cdot \nabla){z}+\nabla \xi+\zeta\nabla\theta=y-y_d,\\
-\nabla\cdot {z}=0,\\
-\frac{\partial \zeta}{\partial t}-\nu_2\Delta\zeta-y\cdot\nabla\zeta={z}\cdot \mathbf{{e}}_2,
\end{gather*}
together with the boundary and initial conditions
\begin{gather*}
\quad {z}=0 ~\text{on}~ \Sigma, \quad \nu_2\frac{\partial\zeta}{\partial \vec{n}}=0 ~\text{on}~ \Sigma.\\
{z}(T)=0,\quad\zeta(T)=0,
\end{gather*}
Then, equality (\ref{sum_equality}) reduces to
$$\iint_{Q}(y-y_d)\cdot \delta y~dxdt=\iint_{\Sigma}\zeta\delta v\chi_{\Sigma_c}~dxdt.$$
It follows from (\ref{def_delta_j}) that
$$\delta J(v)=\iint_{\Sigma_c} (\zeta+\alpha v)\delta v ~ dxdt,$$
which implies that
$$D J(v)=\zeta|_{\Sigma_c}+\alpha v.$$

From the above discussions, we obtain the following results.
\begin{theorem}\label{thm:oc}
	If $u\in L^2(\Sigma_c)$ is an optimal control of  (\ref{OptimalControl})$-$(\ref{objective_functional}), then it satisfies the following optimality system
	\begin{equation}\label{oc}
	DJ(u)=\alpha u+\zeta|_{\Sigma_c}=0,
	\end{equation}
	with $\zeta$ obtained from the solutions of the following fully coupled systems:
	\begin{flalign}\label{state_equation11}
	\left\{
	\begin{aligned}
	&\frac{\partial y}{\partial t}+(y\cdot{\nabla})y-\nu_1 \Delta y+\nabla p=\theta \mathbf{{e}}_2\quad \text{in}\quad \Omega\times(0,T), \\
	&\nabla\cdot y=0\quad \text{in}\quad \Omega\times(0,T),\\
	&y=0\quad \text{on}\quad \Gamma\times(0,T),\\
	&y(0)=y_0,
	\end{aligned}
	\right.
	\end{flalign}
	\begin{flalign}\label{state_equation12}
	\left\{
	\begin{aligned}
	&\frac{\partial \theta}{\partial t}-\nu_2 \Delta\theta+y\cdot\nabla\theta=0\quad \text{in}\quad \Omega\times(0,T), \\
	&\nu_2\frac{\partial\theta}{\partial \vec{n}}=u\chi_{\Sigma_c}\quad \text{on}\quad \Gamma\times(0,T),\\
	&\theta(0)=\theta_0,
	\end{aligned}
	\right.
	\end{flalign}
	and
	\begin{flalign}\label{adjoint_equation1}
	\left\{
	\begin{aligned}
	&-\frac{\partial {z}}{\partial t}-\nu_1\Delta{z}+(\nabla y)^T{z}-(y\cdot \nabla){z}-\nabla \xi+\zeta\nabla\theta=y-y_d\quad \text{in}~\Omega\times(0,T), \\
	&-\nabla\cdot {z}=0\quad\text{in}~\Omega\times(0,T), \\
	&{z}=0\quad \text{on}~ \Gamma\times(0,T),\\
	&{z}(T)=0,
	\end{aligned}
	\right.
	\end{flalign}
	\begin{flalign}\label{adjoint_equation2}
	\left\{
	\begin{aligned}
	&-\frac{\partial \zeta}{\partial t}-\nu_2\Delta\zeta-y\cdot\nabla\zeta={z}\cdot \mathbf{{e}}_2\quad \text{in}\quad \Omega\times(0,T), \\
	&\nu_2\frac{\partial\zeta}{\partial \vec{n}}=0\quad \text{on}\quad \Gamma\times(0,T),\\
	&\zeta(T)=0.\\
	\end{aligned}
	\right.
	\end{flalign}
	
	Here, equations (\ref{adjoint_equation1}) and (\ref{adjoint_equation2}) are the adjoint equations corresponding to velocity and temperature, respectively. Equations (\ref{state_equation11})$-$(\ref{state_equation12}) are the state equations defined in (\ref{state_equation1})$-$(\ref{state_equation2}).
\end{theorem}

\begin{remark}
The optimality systems for the vorticity reduction case (\ref{objective_functional_r}) are the same as those given by (\ref{oc})-(\ref{adjoint_equation2}) but with $y-y_d$ in (\ref{adjoint_equation1}) replaced by $\nabla\times(\nabla\times y)$.
\end{remark}

To compute the gradient $DJ(v)$ for any given $v\in L^2(\Sigma_c)$, it follows from Theorem \ref{thm:oc} that one has to solve the  state system (\ref{state_equation11})-(\ref{state_equation12})  and the adjoint system (\ref{adjoint_equation1})-(\ref{adjoint_equation2}). Both the state and adjoint systems consist of coupled PDEs and they are high-dimensional. Hence, computing the gradient $DJ(v)$ is challenging and computationally expensive, and some sophisticated techniques are required to reduce the computational cost.

\section{Space and time discretizations}\label{se:time_dis}
In this section, we first discuss the $L^2$-projection and the Marchuk-Yanenko splitting for the time discretization of the optimal control problem
(\ref{OptimalControl})--(\ref{objective_functional}). Then, the first-order optimality conditions are derived for the time-discretized optimal control problem.
Finally, a finite element method is presented for the space discretization to obtain a fully discretized formulation of the optimal control problem (\ref{OptimalControl})--(\ref{objective_functional}).

\subsection{An operator splitting method for the time discretization}

Assuming that $T$ is finite, for any given positive integer $N$, let $\Delta t=T/N$ be the time discretization step size and $t_n=n{\Delta t}$,  $\forall n=1,\cdots,N$. Then, we approximate ${U}$ by ${U}^{\Delta t}=[L^2(\Gamma_c)]^{N}$ and $v$ by $\bm{v}=(\bm{v}^n)_{n=1}^N\in {U}^{\Delta t} $.

We can define the scalar product on ${U}^{\Delta t}$ as
$$(\bm{u},\bm{v})_{{U}^{\Delta t}}={\Delta t}\sum_{n=1}^N\int_{\Gamma_c}\bm{u}^n\bm{v}^n dx,\quad \forall \bm{u}, \bm{v} \in {U}^{\Delta t}.$$
%and the functions as
%$$
%\begin{aligned}
%&\begin{aligned}
%T_1(\bm{y}, \bm{z}; \bm{b})&=(\bm{b} \cdot \nabla \bm{y}, \bm{z})+\nu_1 (\nabla \bm{y}, \nabla \bm{z}),
%\end{aligned}\\
%&\begin{aligned}
%T_2(\theta, \zeta; \bm{b})&=(\bm{b} \cdot \nabla \theta, \zeta)+\nu_2 (\nabla \theta, \nabla \zeta).
%\end{aligned}\\
%\end{aligned}
%$$
Then the time-discretized formulation of problem (\ref{OptimalControl})--(\ref{objective_functional}) reads as
\begin{flalign}\label{semi_discrete_optimal_control}
\left\{
\begin{array}{ll}
{\bm{u}}\in {U}^{\Delta t}, \\
J^{\Delta t}({\bm{u}})\leq J^{\Delta t}(\bm{v}),\quad \forall \bm{v}\in {U}^{\Delta t},
\end{array}
\right.
\end{flalign}
with the time-discretized objective functional
\begin{equation*}
J^{\Delta t}(\bm{v})=\frac{1}{2}{\Delta t}\sum_{n=1}^N\int_{\Omega}|\bm{y}^n-y_d^n|^2dx+\frac{\alpha}{2}{\Delta t}\sum_{n=1}^N\int_{\Gamma_c}|\bm{v}^n|^2dx,
\end{equation*}
where the time-discretized target function $y_d^n=y_d(t_n), \forall n=1,\cdots,N$, and $(\bm{y}^n)_{n=1}^N$ are given from $(\bm{v}^n)_{n=1}^N$ by the solution of the following time-discretized state equations:
	\begin{equation}\label{L2StateInitial}
	\bm{y}^0=y_{0},\quad 	\bm{\theta}^0=\theta_0, \quad \bm{p}^0=0,
	\end{equation}
	for $n=1,\cdots, N$,  solve
	\begin{flalign}\label{Advection_Dy}
	\qquad \left\{
	\begin{aligned}
	&\frac{\bm{\tilde{y}}^n-\bm{y}^{n-1}}{\Delta t}-\nu_1 \Delta\bm{\tilde{y}}^n+(\bm{y}^{n-1}\cdot\nabla)\bm{\tilde{y}}^{n}+\nabla \bm{p}^{n-1}= \bm{\theta}^n\bm{e}_2\quad \text{in}\quad \Omega, \\
	&\bm{\tilde{y}}^n=0\quad \text{on}\quad \Gamma,\\
	\end{aligned}
	\right.
	\end{flalign}
	\begin{flalign}\label{Advection_Dtheta}
	\qquad \left\{
	\begin{aligned}
	&\frac{\bm{\theta}^n-\bm{\theta}^{n-1}}{\Delta t}-\nu_2 \Delta\bm{\theta}^n+\bm{y}^{n-1}\cdot\nabla\bm{\theta}^{n}= 0\quad \text{in}\quad \Omega, \\
	&\nu_2\frac{\partial\bm{\theta}^n}{\partial{\vec{n}}}=\bm{v}^n\chi_{\Gamma_c}\quad \text{on}\quad \Gamma,
	\end{aligned}
	\right.
	\end{flalign}
	and
	\begin{flalign}\label{deg_Stokes_y}
	\qquad \left\{
	\begin{aligned}
	&\frac{\bm{y}^{n}-\bm{\tilde{y}}^n}{\Delta t}+\nabla (\bm{p}^n-\bm{p}^{n-1})= 0\quad \text{in}\quad \Omega, \\
	&\nabla\cdot \bm{y}^n=0 \quad \text{in}\quad \Omega, \\
	&\bm{{y}}^n\cdot \vec{n}=0\quad \text{on}\quad \Gamma.
	\end{aligned}
	\right.
	\end{flalign}

%\begin{remark}
%{\color{blue} There is no particular diffculty in applying }
%\end{remark}

Using the above operator splitting scheme (\ref{Advection_Dy})-(\ref{deg_Stokes_y}), we decouple the nonlinearity and the incompressibility in the state equation (\ref{state_equation1})-(\ref{state_equation2}), and meanwhile treat the temperature variable $\theta$ separately. As a result, we only need to handle three rather simple linear equations  (\ref{Advection_Dy})-(\ref{deg_Stokes_y}) at each time step.

\begin{remark}
	Note that the Neumann boundary condition $\bm{y}^ n\cdot\vec{n}=0$ is used to guarantee the well-posedness of (\ref{deg_Stokes_y}).  Alternatively, one may consider replacing (\ref{deg_Stokes_y}) by the following equation with Dirichlet boundary condition:
		\begin{flalign}\label{deg_Stokes_y_2}
		\left\{
		\begin{aligned}
			&\frac{\bm{y}^{n}-\bm{\tilde{y}}^n}{\Delta t}+\nabla (\bm{p}^n-\bm{p}^{n-1})= 0\quad \text{in}\quad \Omega, \\
			&\nabla\cdot \bm{y}^n=0 \quad \text{in}\quad \Omega, \\
			&\bm{{y}}^n =0\quad \text{on}\quad \Gamma.
		\end{aligned}
		\right.
	\end{flalign}
	As mentioned in  \cite[Chapter 7]{glowinski2003finite}, problem \eqref{deg_Stokes_y_2} is generally not well-posed since the boundary condition $\bm{{y}}^n=0$ is too demanding for a solution which does not have the $H^1(\Omega)$-regularity. Despite this fact, the finite element discretized analogue of (\ref{deg_Stokes_y}) (see (\ref{fully-discretized-state-variationalform}))  is well-posed. Moreover, compared with (\ref{Advection_Dy})--(\ref{deg_Stokes_y}), numerical experiments show that more accurate approximate solutions can be always obtained by the scheme (\ref{Advection_Dy}), (\ref{Advection_Dtheta}) and (\ref{deg_Stokes_y_2}).
\end{remark}

\subsection{First-order optimality conditions for (\ref{semi_discrete_optimal_control})--(\ref{deg_Stokes_y}) }\label{se:oc}
In this subsection, we derive the first-order optimality conditions for the time-discretized optimal control problem (\ref{semi_discrete_optimal_control})--(\ref{deg_Stokes_y}) using perturbation analysis as what we have done for the continuous case.

Let $\bm{u}\in{U}^{\Delta t}$ be an optimal control of  (\ref{semi_discrete_optimal_control})-(\ref{deg_Stokes_y}) and $DJ^{\Delta t}(\bm{v})$ be the first-order differential of the functional $J^{\Delta t}$ at $\bm{v}\in\mathcal{U}^{\Delta t}$. Then, the following first-order optimality condition holds
\begin{equation*}
	DJ^{\Delta t}(\bm{u}) = 0.
\end{equation*}
To calculate $DJ^{\Delta t}(\bm{v})$, we first observe that
\begin{equation}\label{delta_J_l2}
	\delta J^{\Delta t}(\bm{v})=(D J^{\Delta t}(\bm{v}),\delta\bm{v})={\Delta t}\sum_{n=1}^N\int_{\Omega}(\bm{y}^n-y_d^n)\cdot\delta\bm{y}^ndx+\alpha{\Delta t}\sum_{n=1}^N\int_{\Gamma_c}\bm{v}^n\delta\bm{v}^ndx,
\end{equation}
and $(\delta\bm{y}^n)_{n=1}^N$, $(\delta\bm{\theta}^n)_{n=1}^N$ and $(\delta\bm{p}^n)_{n=1}^N$ satisfy the perturbed time-discretized state equation:
given
$\delta\bm{y}^0=0,~\delta\bm{\theta}^0=0, ~ \text{and}~	\delta\bm{p}^0=0$, then for $n=1,\cdots, N$,
\begin{gather}
	\left\{
	\begin{aligned}\label{delta_l2_y1}
		&\frac{\delta\bm{\tilde{y}}^n-\delta\bm{y}^{n-1}}{\Delta t}-\nu_1 \Delta\delta\bm{\tilde{y}}^n+(\bm{y}^{n-1}\cdot\nabla)\delta\bm{\tilde{y}}^{n}+(\delta\bm{y}^{n-1}\cdot\nabla)\bm{\tilde{y}}^{n}+ \nabla \delta\bm{p}^{n-1}= \delta\bm{\theta}^n\bm{e}_2\quad \text{in}~ \Omega, \\
		&\delta\bm{\tilde{y}}^n=0, \quad\text{on}~ \Gamma,
	\end{aligned}
	\right.\\
	\left\{
	\begin{aligned}\label{delta_l2_theta}
		&\frac{\delta\bm{\theta}^n-\bm{\delta\theta}^{n-1}}{\Delta t}-\nu_2 \Delta\delta\bm{\theta}^n+\delta\bm{y}^{n-1}\cdot\nabla\bm{\theta}^{n}+\bm{y}^{n-1}\cdot\nabla\delta\bm{\theta}^{n}= 0\quad  \text{in}~\Omega, \\
		&\nu_2\frac{\partial\delta\bm{\theta}^n}{\partial{\vec{n}}}=\delta\bm{v}^n\chi_{\Gamma_c}\text{on}\quad \Gamma.
	\end{aligned}
	\right.\\
	\left\{
	\begin{aligned}\label{delta_l2_y2}
		&\frac{\delta\bm{y}^{n}-\delta\bm{\tilde{y}}^n}{\Delta t}+\nabla( \delta \bm{p}^n- \delta \bm{p}^{n-1} )= 0\quad \text{in}\quad \Omega, \\
		&\nabla\cdot \delta\bm{y}^n=0 \quad \text{in}\quad \Omega, \\
		&\delta\bm{{y}}^n\cdot \vec{n}=0\quad \text{on}\quad \Gamma.
	\end{aligned}
	\right.
\end{gather}

Let us choose variables $(\bm{z}^n)_{n=1}^{N+1}$,$(\tilde{\bm{z}}^n)_{n=1}^{N+1}$, $(\bm{\xi}^n)_{n=1}^N$ and $(\bm{\zeta}^n)_{n=1}^{N+1}$, where $\bm{z}^n,\tilde{\bm{z}}^n, \bm{\xi}^n$ and $\bm{\zeta}^n$ are smooth functions of $x$. Multiplying both sides of the first equation in (\ref{delta_l2_y1}), $(\ref{delta_l2_theta})$ by $\tilde{\bm{z}}^n,$ $\bm{\zeta}^n$, and the first two equations in (\ref{delta_l2_y2}) by $\bm{z}^n$ and $\bm{\xi}^n$, respectively. Then integrating continuously over $\Omega$ and discretely over $(0,T)$, we obtain
\begin{gather}
		\Delta t\sum_{n=1}^N\int_{\Omega}\Big(\frac{\delta\bm{\tilde{y}}^n-\delta\bm{y}^{n-1}}{\Delta t}-\nu_1 \Delta\delta\bm{\tilde{y}}^n+(\bm{y}^{n-1}\cdot\nabla)\delta\bm{\tilde{y}}^{n}+
		(\delta\bm{y}^{n-1}\cdot\nabla)\bm{\tilde{y}}^{n}+\nabla \delta \bm{p}^{n-1} \Big)\cdot  \tilde{\bm{z}}^n dx= \Delta t\sum_{n=1}^N\int_{\Omega}\delta\bm{\theta}^n\bm{e}_2\cdot  \tilde{\bm{z}}^n dx,\label{integ_1}\\
	\Delta t\sum_{n=1}^N\int_{\Omega}\Big(\frac{\delta\bm{\theta}^n-\bm{\delta\theta}^{n-1}}{\Delta t}-\nu_2 \Delta\delta\bm{\theta}^n+\delta\bm{y}^{n-1}\cdot\nabla\bm{\theta}^{n}+\bm{y}^{n-1}\cdot\nabla\delta\bm{\theta}^{n}\Big) \bm{\zeta}^n dx= 0,\label{integ_2}\\
	\Delta t\sum_{n=1}^N\int_{\Omega}\Big(\frac{\delta\bm{y}^{n}-\delta\bm{\tilde{y}}^n}{\Delta t}+\nabla (\delta \bm{p}^n-\delta \bm{p}^{n-1} )\Big)\cdot \bm{z}^n dx= 0,\label{integ_3}\\
	\Delta t\sum_{n=1}^N\int_{\Omega}\nabla  \cdot\delta\bm{y}^n \bm{\xi}^ndx=0.\label{integ_4}
\end{gather}

Since $\delta \bm{y}^0=0$ and $\delta \bm{\theta}^0=0$, it is easy to verify that
\begin{eqnarray}\label{rearrange1}
	\begin{split}
		&\Delta t\sum_{n=1}^N\int_{\Omega}\frac{\delta\bm{\tilde{y}}^n-\delta\bm{y}^{n-1}}{\Delta t}\cdot \tilde{\bm{z}}^n+\frac{\delta\bm{y}^{n}-\delta\bm{\tilde{y}}^n}{\Delta t}\cdot\bm{z}^n dx\\
		=&\Delta t\sum_{n=1}^N\int_{\Omega}\frac{\tilde{\bm{z}}^n-\bm{z}^n}{\Delta t}\delta \tilde{\bm{y}}^n dx+\Delta t\sum_{n=1}^N\int_{\Omega}\frac{{\bm{z}}^n-\tilde{\bm{z}}^{n+1}}{\Delta t}\delta {\bm{y}}^ndx+\int_{\Omega}\tilde{\bm{z}}^{N+1}\delta{\bm{y}}^N dx,
	\end{split}
\end{eqnarray}
and
\begin{equation}\label{rearrange2}
	\Delta t\sum_{n=1}^N\int_{\Omega}\frac{\delta\bm{\theta}^n-\bm{\delta\theta}^{n-1}}{\Delta t}\bm{\zeta}^n dx=\Delta t\sum_{n=1}^N\int_{\Omega}\frac{\bm{\zeta}^n-\bm{\zeta}^{n+1}}{\Delta t}\delta\bm{\theta}^n dx+\int_{\Omega}\bm{\zeta}^{N+1}\delta\bm{\theta}^N dx.
\end{equation}

Using Green's formula and taking the boundary conditions (\ref{delta_l2_y1}), (\ref{delta_l2_theta}) and (\ref{delta_l2_y2}) into account, we have
\begin{eqnarray}\label{rearrange3}
	\begin{split}
		&\Delta t\sum_{n=1}^N\int_{\Omega}	\big((\bm{y}^{n-1}\cdot\nabla)\delta\bm{\tilde{y}}^{n}+(\delta\bm{y}^{n-1}\cdot\nabla)\bm{\tilde{y}}^{n}\big)\cdot \tilde{\bm{z}}^n dx\\
		=&\Delta t\sum_{n=1}^N\int_{\Omega}(\tilde{\bm{y}}^n)^T\tilde{\bm{z}}^n\cdot\delta\bm{y}^{n-1}-(\bm{y}^{n-1}\cdot\nabla)\tilde{\bm{z}}^n\cdot\delta\tilde{\bm{y}}^n -{\nabla \cdot \bm{y}^{n-1} }  \tilde{\bm{z}}^n \cdot\delta\tilde{\bm{y}}^n dx\\
		=&\Delta t\sum_{n=0}^{N-1}\int_{\Omega}(\tilde{\bm{y}}^{n+1})^T\tilde{\bm{z}}^{n+1}\cdot\delta\bm{y}^{n}dx-\Delta t\sum_{n=1}^N\int_{\Omega}((\bm{y}^{n-1}\cdot\nabla)\tilde{\bm{z}}^n+{\nabla \cdot \bm{y}^{n-1} }  \tilde{\bm{z}}^n)\cdot\delta\tilde{\bm{y}}^n dx\\
			=&\Delta t\sum_{n=1}^{N}\int_{\Omega}(\tilde{\bm{y}}^{n+1})^T\tilde{\bm{z}}^{n+1}\cdot\delta\bm{y}^{n}-((\bm{y}^{n-1}\cdot\nabla)\tilde{\bm{z}}^n+{\nabla \cdot \bm{y}^{n-1} }  \tilde{\bm{z}}^n)\cdot\delta\tilde{\bm{y}}^n dx-\Delta t \int_{\Omega}(\tilde{\bm{y}}^{N+1})^T\tilde{\bm{z}}^{N+1}\cdot\delta\bm{y}^{N}dx,
	\end{split}
\end{eqnarray}
where the last equality follows from $\delta \bm{y}^0=0$ and $\tilde{\bm{y}}^{N+1}$ is an auxiliary variable to simplify the derivation of $\tilde{\bm{z}}^{N+1}$.

In a similar way, we can obtain
\begin{eqnarray}\label{rearrange4}
	\begin{split}
		&\Delta t\sum_{n=1}^N\int_{\Omega}\Big(\delta\bm{y}^{n-1}\cdot\nabla\bm{\theta}^{n}+\bm{y}^{n-1}\cdot\nabla\delta\bm{\theta}^{n}\Big) \bm{\zeta}^n dx\\
		=&\Delta t\sum_{n=1}^N\int_{\Omega}\bm{\zeta}^n\nabla\bm{\theta}^{n}\cdot\delta\bm{y}^{n-1}-\bm{y}^{n-1}\cdot\nabla\bm{\zeta}^n\delta\bm{\theta}^{n} -({\nabla \cdot \bm{y}^{n-1} })\bm{\zeta}^n\delta\bm{\theta}^{n} dx\\
		= &\Delta t\sum_{n=1}^{N-1}\int_{\Omega}\bm{\zeta}^{n+1}\nabla\bm{\theta}^{n+1}\cdot\delta\bm{y}^{n}dx- \Delta t\sum_{n=1}^N\int_{\Omega}(\bm{y}^{n-1}\cdot\nabla\bm{\zeta}^n+{\nabla \cdot \bm{y}^{n-1} }\bm{\zeta}^n)\delta\bm{\theta}^{n} dx\\
		= &\Delta t\sum_{n=1}^{N}\int_{\Omega}\bm{\zeta}^{n+1}\nabla\bm{\theta}^{n+1}\cdot\delta\bm{y}^{n}-(\bm{y}^{n-1}\cdot\nabla\bm{\zeta}^n+{\nabla \cdot \bm{y}^{n-1} }\bm{\zeta}^n)\delta\bm{\theta}^{n} dx-\int_{\Omega}\bm{\zeta}^{N+1}\nabla\bm{\theta}^{N+1}\cdot\delta\bm{y}^{N}dx \\
	\end{split}
\end{eqnarray}
and
\begin{eqnarray}\label{rearrange5}
	\begin{split}
		&\Delta t\sum_{n=1}^N\int_{\Omega}\nabla \delta\bm{p}^{n-1}\cdot\tilde{\bm{z}}^n+\nabla \delta(\bm{p}^{n}-\bm{p}^{n-1})\cdot {\bm{z}}^n dx\\
		=&\Delta t\sum_{n=0}^{N-1}\int_{\Omega}\nabla \cdot(\bm{z}^n+\bm{\tilde{z}}^{n+1}-\bm{z}^{n+1}) \delta \bm{p}^{n} dx+\int_{\Omega}\nabla \cdot \bm{z}^N \delta \bm{p}^{N} dx\\
		&-\Delta t\sum_{n=0}^{N-1}\int_{\Gamma}\vec{n} \cdot(\bm{z}^n+\bm{\tilde{z}}^{n+1}-\bm{z}^{n+1}) \delta \bm{p}^{n} dx-\int_{\Gamma}\vec{n} \cdot \bm{z}^N \delta \bm{p}^{N} dx\\
	=&\Delta t\sum_{n=1}^{N}\int_{\Omega}\nabla \cdot(\bm{z}^n+\bm{\tilde{z}}^{n+1}-\bm{z}^{n+1}) \delta \bm{p}^{n} dx-\int_{\Omega}\nabla \cdot (\bm{\tilde{z}}^{N+1}-\bm{z}^{N+1})\delta \bm{p}^{N} dx\\
	&-\Delta t\sum_{n=1}^{N}\int_{\Gamma}\vec{n}  \cdot(\bm{z}^n+\bm{\tilde{z}}^{n+1}-\bm{z}^{n+1}) \delta \bm{p}^{n} dx+\int_{\Gamma}\vec{n}  \cdot (\bm{\tilde{z}}^{N+1}-\bm{z}^{N+1})\delta \bm{p}^{N} dx.
	\end{split}
\end{eqnarray}
Substituting (\ref{rearrange1})$-$(\ref{rearrange5}) into (\ref{integ_1})$-$(\ref{integ_4}) and summing them up, and then applying Green's formula to the rest of terms, we obtain
\begin{eqnarray}\label{sum_three_eqn}
	\begin{aligned}
	&\Delta t\sum_{n=1}^{N}\int_{\Omega} \Big(\frac{{\bm{z}}^n-\tilde{\bm{z}}^{n+1}}{\Delta t}\cdot\delta {\bm{y}}^n-\nabla\bm{\xi}^n\cdot\delta\bm{y}^n+(\tilde{\bm{y}}^{n+1})^T\tilde{\bm{z}}^{n+1}\cdot\delta\bm{y}^{n}+\bm{\zeta}^{n+1}\nabla\bm{\theta}^{n+1}\cdot\delta\bm{y}^{n}\\
	&+\frac{\tilde{\bm{z}}^n-\bm{z}^n}{\Delta t}\cdot\delta \tilde{\bm{y}}^n-\nu_1\Delta\tilde{\bm{z}}^n\cdot\delta \tilde{\bm{y}}^n -((\bm{y}^{n-1}\cdot\nabla)\tilde{\bm{z}}^n+{\nabla \cdot \bm{y}^{n-1} }  \tilde{\bm{z}}^n)\cdot\delta\tilde{\bm{y}}^n\\
	&+\frac{\bm{\zeta}^n-\bm{\zeta}^{n+1}}{\Delta t}\delta\bm{\theta}^n-\nu_2\Delta\bm{\zeta}^n\delta{\bm{\theta}}^n-(\bm{y}^{n-1}\cdot\nabla\bm{\zeta}^n+{\nabla \cdot \bm{y}^{n-1} }\bm{\zeta}^n)\delta\bm{\theta}^{n} +\nabla \cdot(\bm{z}^n+\bm{\tilde{z}}^{n+1}-\bm{z}^{n+1}) \delta \bm{p}^{n} \Big)dx\\
		&+\Delta t \int_{\Omega}\Big(\tilde{\bm{z}}^{N+1}\cdot\delta{\bm{y}}^N+\bm{\zeta}^{N+1}\delta\bm{\theta}^N-(\tilde{\bm{y}}^{N+1})^T\tilde{\bm{z}}^{N+1}\cdot\delta\bm{y}^{N}-\bm{\zeta}^{N+1}\nabla\bm{\theta}^{N+1}\cdot\delta\bm{y}^{N}-\nabla \cdot (\bm{\tilde{z}}^{N+1}-\bm{z}^{N+1})\delta \bm{p}^{N} \Big)dx\\
		=&  \Delta t\sum_{n=1}^N\Big(\int_{\Omega}\tilde{\bm{z}}^n \cdot \bm{e}_2  \delta\bm{\theta}^ndx+\nu_1\int_{\Gamma}\frac{\partial\delta \tilde{\bm{y}}^n}{\partial \vec{n}} \tilde{\bm{z}}^n-\frac{\partial \tilde{\bm{z}}^{n} }{\partial \vec{n}}\cdot\delta \tilde{\bm{y}}^n dx+\nu_2\int_{\Gamma}\frac{\partial\delta\bm{\theta}^n}{\partial \vec{n}}\bm{\zeta}^n-\frac{\partial \bm{\zeta}^n}{\partial \vec{n}}\cdot\delta \bm{\theta}^n dx\Big)\\
		&+\Delta t\sum_{n=1}^{N}\int_{\Omega}\vec{n}  \cdot(\bm{z}^n+\bm{\tilde{z}}^{n+1}-\bm{z}^{n+1}) \delta \bm{p}^{n} dx-\int_{\Omega}\vec{n}  \cdot (\bm{\tilde{z}}^{N+1}-\bm{z}^{N+1})\delta \bm{p}^{N} dx.
	\end{aligned}
\end{eqnarray}
Let $(\tilde{\bm{z}}^n)_{n=1}^{N+1}$,  $({\bm{z}}^n)_{n=1}^{N+1}$ and $(\bm{\zeta}^n)_{n=1}^{N+1}$ satisfy the following equations
\begin{equation}\label{L2adjointInitial}
	\tilde{\bm{z}}^{N+1}=0,\quad {\bm{z}}^{N+1}=0  \quad\text{and}\quad \bm{\zeta}^{N+1}=0,
\end{equation}
for $n=N,N-1,\cdots,2,1$
\begin{flalign}\label{deg_Stokes_z}
	\left\{
	\begin{aligned}
		&\frac{{\bm{z}}^n-\tilde{\bm{z}}^{n+1}}{\Delta t}-\nabla\bm{\xi}^n+(\tilde{\bm{y}}^{n+1})^T\tilde{\bm{z}}^{n+1}+\bm{\zeta}^{n+1}\nabla\bm{\theta}^{n+1}=\bm{y}^n-\bm{y}_d^n\quad \text{in}\quad \Omega, \\
		&\nabla\cdot(\bm{z}^n+\bm{\tilde{z}}^{n+1}-\bm{z}^{n+1})=0\quad \text{in}\quad \Omega, \\
		&\bm{z}^n\cdot \vec{n}=0\quad \text{on}\quad \Gamma,
	\end{aligned}
	\right.
\end{flalign}
\begin{flalign}\label{Advection_Dz}
	\left\{
	\begin{aligned}
		&\frac{\tilde{\bm{z}}^n-\bm{z}^n}{\Delta t}-\nu_1\Delta\tilde{\bm{z}}^n-((\bm{y}^{n-1}\cdot\nabla)\tilde{\bm{z}}^n+{\nabla \cdot \bm{y}^{n-1} }  \tilde{\bm{z}}^n)=0\quad \text{in}\quad \Omega, \\
		&\tilde{\bm{z}}^n=0\quad \text{on}\quad \Gamma
	\end{aligned}
	\right.
\end{flalign}
and
\begin{flalign}\label{Advection_Dzeta}
	\left\{
	\begin{aligned}
		&\frac{\bm{\zeta}^n-\bm{\zeta}^{n+1}}{\Delta t}-\nu_2\Delta\bm{\zeta}^n-(\bm{y}^{n-1}\cdot\nabla\bm{\zeta}^n+{\nabla \cdot \bm{y}^{n-1} }\bm{\zeta}^n)=\tilde{\bm{z}}^n \cdot \bm{e}_2\quad \text{in}\quad \Omega, \\
		&\nu_2\frac{\partial\bm{\zeta}^n}{\partial \vec{n}}=0\quad \text{on}\quad \Gamma.
	\end{aligned}
	\right.
\end{flalign}
Then, equality (\ref{sum_three_eqn}) reduces to
\begin{equation*}
	\Delta t\sum_{n=1}^{N}\int_{\Omega}(\bm{y}^n-\bm{y}_d^n)\delta\bm{y}^ndx=
	\Delta t\sum_{n=1}^{N}\int_{\Gamma}\bm{\zeta}^n\delta\bm{v}^n\chi_{\Gamma_c}dx.
\end{equation*}
Taking (\ref{delta_J_l2}) into account, we obtain
\begin{equation*}
	\delta J^{\Delta t}(\bm{v})=(D J^{\Delta t}(\bm{v}),\delta\bm{v})=\Delta t\sum_{n=1}^{N}\int_{\Gamma_c}(\bm{\zeta}^n+\alpha\bm{v}^n)\delta\bm{v}^ndx,
\end{equation*}
which implies that
\begin{equation}\label{grad}
	DJ^{\Delta t}(\bm{v})=(\bm{\zeta}^n\big|_{\Gamma_c}+\alpha\bm{v}^n)_{n=1}^N.
\end{equation}

We summarize the first-order optimality conditions for the time-discretized problem (\ref{semi_discrete_optimal_control})-(\ref{deg_Stokes_y}) as follows:
\begin{theorem}
	If $\bm{u}=\{\bm{u}^n\}_{n=1}^N\in{U}^{\Delta t}$ is an optimal control of the time-discretized problem (\ref{semi_discrete_optimal_control})-(\ref{deg_Stokes_y}), then the following optimality conditions hold
	\begin{equation*}
		\alpha\bm{u}^n+\bm{\zeta}^n|_{\Gamma_c}=0, \quad\forall n=1,2,\cdots,N,
	\end{equation*}
	where $(\bm{\zeta}^n)_{n=1}^N$ is the successive solution of the time-discretized state equations (\ref{L2StateInitial})$-$(\ref{deg_Stokes_y}) and the corresponding adjoint equations (\ref{L2adjointInitial})$-$(\ref{Advection_Dzeta}).
\end{theorem}

To compute a gradient $DJ^{\Delta t}(\bm{v})$, we only need to solve four linear advection-diffusion equations: (\ref{Advection_Dy}), (\ref{Advection_Dtheta}), (\ref{Advection_Dz})  and (\ref{Advection_Dzeta}), and two degenerated Stokes equations: (\ref{deg_Stokes_y}) and (\ref{deg_Stokes_z}), at each time step.  All these equations can be easily solved by some well-developed numerical methods in the literature.  For instance,  a fixed-point iterative process is proposed in \cite{BN99} for solving advection-diffusion problems,  and a preconditioned conjugate gradient method is suggested in \cite{glowinski2003finite} to solve degenerated Stokes equations. More details can be found in \cite{glowinski2003finite} and references therein.

{
\begin{remark}\label{br1}
	Following \cite{glowinski2003finite}, implementing the Peaceman-Rachford splitting method \cite{peaceman1955} to (\ref{state_equation1})-(\ref{state_equation2}) yields the following time-discretized equations (with $ 0 < \beta < 1$, $0 < \gamma < 1$ and
	$\beta + \gamma = 1$): $\bm{y}^0=y_{0}, \bm{\theta}^0=\theta_0, $ then for $n=1,\cdots, N$,  solve
	\begin{flalign}\label{prsm1}
		\qquad \left\{
		\begin{aligned}
			&\frac{\bm{\tilde{y}}^n-\bm{y}^{n-1}}{\Delta t/2}-\beta\nu_1 \Delta\bm{\tilde{y}}^n+\nabla \bm{\tilde{p}}^{n}= \bm{\theta}^{n-1}\bm{e}_2+\gamma\nu_1 \Delta\bm{y}^{n-1}-(\bm{y}^{n-1}\cdot\nabla)\bm{{y}}^{n-1}\quad \text{in}\quad \Omega, \\
			&\nabla\cdot \bm{\tilde{y}}^n=0 \quad \text{in}\quad \Omega, \\
			&\bm{\tilde{y}}^n=0\quad \text{on}\quad \Gamma,\\
		\end{aligned}
		\right.
	\end{flalign}
	\begin{flalign}\label{prsm2}
	\qquad \left\{
	\begin{aligned}
		&\frac{\bm{y}^{n}-\bm{\tilde{y}}^n}{\Delta t/2}-\gamma\nu_1 \Delta\bm{{y}}^n+(\bm{y}^{n}\cdot\nabla)\bm{{y}}^{n}= \bm{\theta}^{n-1}\bm{e}_2+\beta\nu_1 \Delta\bm{\tilde{y}}^n-\nabla\bm{\tilde{p}}^n\quad \text{in}\quad \Omega, \\
		&\nabla\cdot \bm{y}^n=0 \quad \text{in}\quad \Omega, \\
		&\bm{{y}}^n\cdot \vec{n}=0\quad \text{on}\quad \Gamma,
	\end{aligned}
	\right.
\end{flalign}
and
	\begin{flalign}\label{prsm3}
		\qquad \left\{
		\begin{aligned}
			&\frac{\bm{\theta}^n-\bm{\theta}^{n-1}}{\Delta t}-\nu_2 \Delta\bm{\theta}^n+\bm{y}^{n}\cdot\nabla\bm{\theta}^{n}= 0\quad \text{in}\quad \Omega, \\
			&\nu_2\frac{\partial\bm{\theta}^n}{\partial{\vec{n}}}=\bm{v}^n\chi_{\Gamma_c}\quad \text{on}\quad \Gamma.
		\end{aligned}
		\right.
	\end{flalign}
We can see that solving the resulting subproblems (\ref{prsm1})-(\ref{prsm3}) is not difficult.  However, in a similar way as what we have done for deriving (\ref{grad}), we found that it is challenging to obtain the adjoint equations associated with the time-discretized equations (\ref{prsm1})-(\ref{prsm3}).  A particular reason is that some terms arise in different time intervals (e.g.,  $\beta\nu_1 \Delta\bm{\tilde{y}}^n$ and $\nabla \bm{\tilde{p}}^{n}$ arise in both $[t^{n-1}, t^{n-1}+\Delta t/2]$ and $[ t^{n-1}+\Delta t/2, t^n]$),  which makes their adjoint variables difficult to be determined. The same concerns apply also to the Douglas-Rachford splitting method \cite{douglas1956}  and the $\theta$-scheme \cite{glowinski1985,glowinski1986}. By contrast, there is no coupled term between different time intervals in our proposed scheme (\ref{Advection_Dy})-(\ref{deg_Stokes_y}), and thus no particular difficulty is introduced in deriving the corresponding adjoint equation as shown in (\ref{delta_J_l2})-(\ref{grad}).
\end{remark}
}

Using the notations introduced in Section \ref{se:notations}, it is easy to show that the variational formulations of the time-discretized state equations (\ref{L2StateInitial})$-$(\ref{deg_Stokes_y}) read as, for $n=1,2,\cdots, N$, find $(\bm{y}^n,\tilde{\bm{y}}^n,\bm{p}^n,\bm{\theta}^n)\in [H^1(\Omega)]^2\times [H_0^1(\Omega)]^2\times P\times H^1(\Omega)$, such that
\begin{flalign*}%\label{time-discretized-state-variationalform}
	\left\{
	\begin{aligned}
		&(\frac{\tilde{\bm{y}}^n-\bm{y}^{n-1}}{\Delta t}, \varphi)-b(\bm{p}^{n-1},{\varphi})+  a(\tilde{\bm{y}}^n, \varphi)+{c}(\bm{y}^{n-1}, \tilde{\bm{y}}^n, \varphi)=(\bm{\theta}^n \bm{e}_2, \varphi),\quad \forall \varphi\in [H_0^1(\Omega)]^2 ,\\
		&(\frac{\bm{\theta}^n-\bm{\theta}^{n-1}}{\Delta t},\psi)+d(\bm{\theta}^n,\psi)+{e}(\bm{y}^{n-1},\bm{\theta}^n, \psi)=(\bm{v}^{n}\chi_{\Gamma_c},\psi)_\Gamma,\quad \forall \psi\in H^1(\Omega),\\
		&(\frac{\bm{y}^{n}-\tilde{\bm{y}}^n}{\Delta t}, {\varphi})-b(\bm{p}^n-\bm{p}^{n-1},{\varphi})=0,\quad \forall {\varphi}\in [H^1(\Omega)]^2,\\
		&b(q,\bm{y}^{n})=0,\quad \forall q\in P,\\
		&\bm{y}^{0}=y_0,\quad \bm{\theta}^{0}=\theta_0.
	\end{aligned}
	\right.
\end{flalign*}

Similarly, the variational formulations of the adjoint equations (\ref{L2adjointInitial})$-$(\ref{Advection_Dzeta}) are then given by: for $n=N,N-1,\cdots,1$, find $(\bm{z}^n,\tilde{\bm{z}}^n,\bm{\xi}^n,\bm{\zeta}^n)\in [H^1(\Omega)]^2\times [H_0^1(\Omega)]^2\times P\times H^1(\Omega)$, such that
\begin{flalign*}%\label{time-discretized-adjoint-variationalform}
	\left\{
	\begin{aligned}
		&(\frac{\bm{z}^{n}-\tilde{\bm{z}}^{n+1}}{\Delta t}, \varphi)+b(\bm{\xi}^n,\varphi)+{c}(\varphi,\tilde{\bm{y}}^{n+1}, \tilde{\bm{z}}^{n+1})\\
		&\qquad\qquad\qquad\quad+{e}(\varphi,\bm{\theta}^{n+1},\bm{\zeta}^{n+1})=(\bm{y}^n-y_d^n,\varphi),\forall \varphi\in [H^1(\Omega)]^2,\\
		&-b(q,\bm{z}^n+\bm{\tilde{z}}^{n+1}-\bm{z}^{n+1})=0, \forall q\in P,\\
		&(\frac{\tilde{\bm{z}}^n-\bm{z}^{n}}{\Delta t}, \varphi)+a(\tilde{\bm{z}}^n,\varphi)+{c}(\bm{y}^{n-1},\varphi,\tilde{\bm{z}}^n)=0, \forall \varphi\in [H^1_0(\Omega)]^2,\\
		&(\frac{\bm{\zeta}^n-\bm{\zeta}^{n+1}}{\Delta t},\psi)+d(\bm{\zeta}^n,\psi)+{e}(\bm{y}^{n-1},\psi,\bm{\zeta}^n )=(\tilde{\bm{z}}^n\cdot \bm{e}_2,\psi), \forall \psi\in H^1(\Omega),\\
		&\tilde{\bm{z}}^{N+1}=0, \quad {\bm{z}}^{N+1}=0 ,\quad \bm{\zeta}^{N+1}=0.
	\end{aligned}
	\right.
\end{flalign*}

\subsection{Space discretization}
	   For the space discretization of  problem (\ref{semi_discrete_optimal_control})-(\ref{deg_Stokes_y}), we employ the Bercovier-Pironneau finite element pair \cite{BP79} (a.k.a $P_1$-$P_1$ iso $P_2$ finite element).  More concretely, the velocity $y$ and the pressure $p$ are approximated by conforming linear finite elements with the mesh sizes $h$ and $2h$, respectively.  For the approximation of the temperature $\theta$ and the control $v$, we use the same finite element space as the one used for the velocity $y$.
	
	   For simplicity, we suppose from now on that $\Omega$ is a polygonal
	   domain of $\mathbb{R}^2$ (or has been approximated by a
	   family of such domains). Let $\mathcal{T}_H$ be a classical triangulation of $\Omega$, with $H$ the largest length of the edges of the triangles of $\mathcal{T}_H$. From $\mathcal{T}_{H}$ we construct $\mathcal{T}_{h}$ with $h=H/2$ by joining the mid-points of the
	   edges of the triangles of $\mathcal{T}_{H}$.
	   We consider three finite element spaces $V_h$, $V_{0h}$ and $P_H$ defined by
	   \begin{eqnarray*}
	   	&&V_h = \{\varphi_h| \varphi_h\in C^0(\bar{\Omega}); \varphi_h{\mid}_{\mathbb{T}}\in \mathcal{P}_1, \forall\, {\mathbb{T}}\in\mathcal{T}_h\},\\
   	&&V_{0h} =\{\varphi_h| \varphi_h\in V_h, \varphi_h{\mid}_{\Gamma}=0\}:=V_h\cap H_0^1(\Omega),\\
   	&& P_H=\{q_H|q_H \in C^0(\bar{\Omega}); q_H{\mid}_{\mathbb{T}}\in \mathcal{P}_1, \forall\, {\mathbb{T}}\in\mathcal{T}_H;\int_{\Omega}q_Hdx=0\},
   \end{eqnarray*}
	   with $\mathcal{P}_1$ the space of the polynomials of two variables of degree $\leq 1$.

	   With the above finite element spaces, the fully discretized formulation of problem (\ref{OptimalControl})-(\ref{objective_functional}) reads as
	   \begin{flalign}\label{fully_discrete_optimal_control}
	   	\left\{
	   	\begin{array}{ll}
	   		{\bm{u}_h}\in [V_h]^N, \\
	   		J_h^{\Delta t}({\bm{u}_h})\leq J_h^{\Delta t}(\bm{v}_h),\quad \forall \bm{v}_h\in [V_h]^N,
	   	\end{array}
	   	\right.
	   \end{flalign}
	   with the time-discretized objective functional
	   \begin{equation}\label{fully_discrete_L2Objective}
	   	J_h^{\Delta t}(\bm{v}_h)=\frac{1}{2}{\Delta t}\sum_{n=1}^N\int_{\Omega}|\bm{y}_h^n-y_{dh}^n|^2dx+\frac{\alpha}{2}{\Delta t}\sum_{n=1}^N\int_{\Gamma_c}|\bm{v}_h^n|^2dx,
	   \end{equation}
	   where $(\bm{y}_h^n)_{n=1}^N$ are given from $(\bm{v}^n_h)_{n=1}^N$ by the solution of the following fully discretized state equations:
	   for $n=1,2,\cdots,N$, find $(\bm{y}_h^n,\tilde{\bm{y}}_h^n,\bm{p}_h^n,\bm{\theta}_h^n)\in [V_h]^2\times [V_{0h}]^2\times P_H\times V_h$, such that
	   \begin{equation}\label{fully-discretized-state-variationalform}
	\left\{
	\begin{aligned}
		&(\frac{\tilde{\bm{y}}_h^n-\bm{y}_h^{n-1}}{\Delta t}, \varphi_h)-b(\bm{p}_H^{n-1},{\varphi_h})+  a(\tilde{\bm{y}}_h^n, \varphi_h)+{c}(\bm{y}_h^{n-1}, \tilde{\bm{y}}_h^n, \varphi_h)=(\bm{\theta}_h^n \bm{e}_2, \varphi_h),\quad \forall \varphi_h\in [V_{0h}]^2 ,\\
		&(\frac{\bm{\theta}_h^n-\bm{\theta}_h^{n-1}}{\Delta t},\psi_h)+d(\bm{\theta}_h^n,\psi)+{e}(\bm{y}_h^{n-1},\bm{\theta}_h^n, \psi_h)=(\bm{v}_h^{n}\chi_{\Gamma_c},\psi_h)_\Gamma,\quad \forall \psi_h\in V_h,\\
		&(\frac{\bm{y}_h^{n}-\tilde{\bm{y}}_h^n}{\Delta t}, {\varphi_h})-b(\bm{p}_H^n-\bm{p}_H^{n-1},{\varphi_h})=0,\quad \forall {\varphi_h}\in [V_h]^2,\\
		&b(q_H,\bm{y}_h^{n})=0,\quad \forall q_H\in P_H,\\
		&\bm{y}_h^{0}=y_{0h},\quad \bm{\theta}_h^{0}=\theta_{0h}.
	\end{aligned}
	\right.
      \end{equation}

\begin{remark}
	 The fully discretized state equations corresponding to the scheme (\ref{Advection_Dy}), (\ref{Advection_Dtheta}) and (\ref{deg_Stokes_y_2}) read as:
	for $n=1,2,\cdots,N$, find $(\bm{y}_h^n,\tilde{\bm{y}}_h^n,\bm{p}_h^n,\bm{\theta}_h^n)\in [V_{0h}]^2\times [V_{0h}]^2\times P_H\times V_h$, such that
	\begin{flalign}\label{fully-discretized-state-variationalform0}
		\left\{
		\begin{aligned}
			&(\frac{\tilde{\bm{y}}_h^n-\bm{y}_h^{n-1}}{\Delta t}, \varphi_h)-b(\bm{p}_H^{n-1},{\varphi_h})+  a(\tilde{\bm{y}}_h^n, \varphi_h)+{c}(\bm{y}_h^{n-1}, \tilde{\bm{y}}_h^n, \varphi_h)=(\bm{\theta}_h^n \bm{e}_2, \varphi_h),\quad \forall \varphi_h\in [V_{0h}]^2 ,\\
			&(\frac{\bm{\theta}_h^n-\bm{\theta}_h^{n-1}}{\Delta t},\psi_h)+d(\bm{\theta}_h^n,\psi)+{e}(\bm{y}_h^{n-1},\bm{\theta}_h^n, \psi_h)=(\bm{v}_h^{n}\chi_{\Gamma_c},\psi_h)_\Gamma,\quad \forall \psi_h\in V_h,\\
			&(\frac{\bm{y}_h^{n}-\tilde{\bm{y}}_h^n}{\Delta t}, {\varphi_h})-b(\bm{p}_H^n-\bm{p}_H^{n-1},{\varphi_h})=0,\quad \forall {\varphi_h}\in [V_{0h}]^2,\\
			&b(q_H,\bm{y}_h^{n})=0,\quad \forall q_H\in P_H,\\
			&\bm{y}_h^{0}=y_{0h},\quad \bm{\theta}_h^{0}=\theta_{0h}.
		\end{aligned}
		\right.
	\end{flalign}
\end{remark}

 	In a similar way as what we have done in Section \ref{se:time_dis},  it is easy to derive the first-order optimality conditions for the fully discretized optimal control problem (\ref{fully_discrete_optimal_control})-(\ref{fully_discrete_L2Objective}) constrained by (\ref{fully-discretized-state-variationalform}) or  (\ref{fully-discretized-state-variationalform0}). For succinctness, we omit the details. 	 \

\section{A BFGS method for problem (\ref{fully_discrete_optimal_control})-(\ref{fully-discretized-state-variationalform})}
\label{sec:Numerical algorithm}
Quasi-Newton methods are milestones in the development of numerical optimization, and they are particularly efficient for large-scale problems. The common idea of quasi-Newton methods is approximating the Hessian or its inverse with only gradient information. Thanks to its nice properties and excellent performance, the BFGS method is the most representative and widely-used quasi-Newton algorithm, see, e.g. \cite{nocedal2006sequential}. 
In this section, we shall discuss the application of the BFGS method to the solution of the fully discretized optimal control problem (\ref{fully_discrete_optimal_control})-(\ref{fully-discretized-state-variationalform}), and propose an easily implementable algorithm tailored to the structure of (\ref{fully_discrete_optimal_control})-(\ref{fully-discretized-state-variationalform}). For ease of notations, the subscripts $h$ in all variables are dropped.

\subsection{A generic BFGS method}\label{se:generic_bfgs}
Conceptually, implementing the BFGS to problem (\ref{fully_discrete_optimal_control})-(\ref{fully-discretized-state-variationalform}), we readily obtain Algorithm \ref{alg:CBFGS}.
\begin{algorithm}
	\caption{A generic BFGS method for the solution of problem (\ref{fully_discrete_optimal_control})-(\ref{fully-discretized-state-variationalform})}
	\label{alg:CBFGS}
\begin{algorithmic}[1]
		 \STATE Given an initial point $\bm{u}_0$, an initial inverse Hessian approximation $\mathcal{H}_0$.
		\STATE Compute $\bm{g}_0=D J_h^{\Delta t}(\bm{u}_0)$; If $\bm{g}_0=0$, then $\bm{u}=\bm{u}_0$, else
        \WHILE {not converged}
        \STATE Compute the search direction
        $\bm{d}_k=-\mathcal{H}_k\bm{g}_k$ and a step size $\rho_k$.
        \STATE Update $\bm{u}_{k+1}=\bm{u}_k+\rho_k \bm{d}_k$ and $\bm{g}_{k+1}=D J_h^{\Delta t}(\bm{u}_{k+1})$.
        \STATE Compute $\delta \bm{u}_k=\bm{u}_{k+1}-\bm{u}_k$ and $\delta \bm{g}_k=\bm{g}_{k+1}-\bm{g}_k$.
        \STATE Update $\mathcal{H}_{k+1}$ via
      $$\mathcal{H}_{k+1}=
        V_k^T\mathcal{H}_kV_k+\eta_k \delta \bm{u}_k\delta \bm{u}_k^T,$$
        where $V_k=I-\eta_k \delta \bm{g}_k  \delta \bm{u}_k^T,$ and $ \eta_k=\frac{1}{\delta\bm{g}_k^T\delta\bm{u}_k}.$
        \ENDWHILE
 \end{algorithmic}
\end{algorithm}

The direct application of the BFGS method to problem (\ref{fully_discrete_optimal_control})-(\ref{fully-discretized-state-variationalform}), however, is not practically implementable. In particular, it is necessary to specify the strategies for determining an appropriate step size $\rho_k$ and a practical inverse  Hessian approximation $\mathcal{H}_k$ at each iteration. Recall that the discretized problem  (\ref{fully_discrete_optimal_control})-(\ref{fully-discretized-state-variationalform}) is high-dimensional and the resulting Hessian is usually dense. It is thus necessary to find efficient strategies for determining the step size $\rho_k$ and the inverse Hessian approximation $\mathcal{H}_k$ to implement the BFGS method or its variants. Because of the huge dimensionality and thus the demanding requirement of memory in computation, we focus on the classic limited-memory BFGS (L-BFGS) method in \cite{nocedal1980} and elaborate on these issues in the following part of this section.

\subsection{Computation of the step size $\rho_k$}\label{se:bfgs_step}

A crucial step to implement Algorithm \ref{alg:CBFGS} is the computation of an appropriate step size $\rho_k$. That is, finding $\rho$ to substantially reduce the univariate function
$$
P_k(\rho):=J_h^{\Delta t}(\bm{u}_k+\rho \bm{d}_k),
$$
where $\bm{u}_k, \bm{d}_k \in [V_h]^N$ are known.
To this end,  a natural idea is to employ some line search strategies, such as the backtracking strategy based on the Armijo--Goldstein condition or the Wolf condition, see e.g., \cite{nocedal2006sequential}. It is worth noting that these line search strategies require the evaluation of $J_h^{\Delta t}(\bm{v})$ repeatedly, which is computationally expensive because every evaluation of $J_h^{\Delta t}(\bm{v})$ for a given $\bm{v}\in [V_h]^N$ requires solving the state equation (\ref{fully-discretized-state-variationalform}).  To address this issue, we advocate the following step size seeking strategy, and a similar idea can also be found in \cite{GSYY2022}.

We consider linearizing the fully discretized state equation (\ref{fully-discretized-state-variationalform}) to find an appropriate step size.
Recall that in the  fully discretized  optimal control problem (\ref{fully_discrete_optimal_control})-(\ref{fully-discretized-state-variationalform}), the states $(\bm{y},\bm{\tilde{y}},\bm{\theta})=(\bm{y}^n,\bm{\tilde{y}}^n,\bm{\theta}^n)_{n=1}^N$ are given from $\bm{v}=(\bm{v}^n)_{n=1}^N$ by the solution of the fully discretized state equation (\ref{fully-discretized-state-variationalform}). We introduce the control-to-sate operator $S_h^{\Delta t}$ associated with the equation (\ref{fully-discretized-state-variationalform}), which maps $\bm{v}$ to $\bm{y}=S_h^{\Delta t}(\bm{v})$. Then the first-order approximation of mapping $\rho\rightarrow S_h^{\Delta t}(\bm{u}_k+\rho \bm{d}_k)$  at $\rho=0$ is 
$$S_h^{\Delta t}(\bm{u}_k)+\rho {S_h^{\Delta t}}'(\bm{u}_k)\bm{d}_k:=\bm{y}_k+\rho \bm{w}_k,$$ 
where $\bm{w}_k$
 are given from $\bm{u}_k$, $\bm{d}_k$ and $(\bm{y}_k,\bm{\tilde{y}}_k,\bm{\theta}_k)$ by the solutions of the following equations:
	\begin{flalign*}%\label{Advection_Dy_linearized}
    \left\{
	\begin{aligned}
	&	\bm{w}_k^0=0,\quad \bm{\phi}_k^0=0 \text{~and~} \bm{q}_{k}^0=0,\\
	&\int_{\Omega}\frac{\bm{\tilde{w}}_k^n-\bm{w}_k^{n-1}}{\Delta t}\varphi_h dx+\int_{\Omega}\nu_1 \nabla\bm{\tilde{w}}_k^n\cdot\nabla \varphi_hdx\\
	&\qquad\qquad+\int_{\Omega}\left((\bm{y}_k^{n-1}\cdot\nabla)\bm{\tilde{w}}_k^{n}+(\bm{w}_k^{n-1}\cdot\nabla)\bm{\tilde{y}}_k^{n}+\nabla \bm{q}_{k}^{n-1}\right)\varphi_hdx= \int_{\Omega}\bm{\phi}_k^n\bm{e}_2\cdot \varphi_hdx, \forall \varphi_h\in [V_{0h}]^2,\\
	&\int_{\Omega} \frac{\bm{\phi}_k^n-\bm{\phi}_k^{n-1}}{\Delta t}\psi_hdx+\nu_2\int_{\Omega} \nabla\bm{\phi}_k^n\cdot\nabla \psi_hdx+\int_{\Omega}\left(\bm{w}_k^{n-1}\cdot\nabla \bm{\theta}_k^n+\bm{y}_k^{n-1}\cdot\nabla\bm{\phi}_k^{n}\right)\psi_hdx= \int_{\Gamma_c}\bm{d}_k^n\psi_hdx, \forall \psi_h\in V_{h},\\
	&\int_{\Omega}\frac{\bm{w}_k^{n}-\bm{\tilde{w}}_k^n}{\Delta t}\varphi_hdx+\int_{\Omega}\nabla (\bm{q}_k^n-\bm{q}_{k}^{n-1})\varphi_hdx= 0, \forall \varphi_h\in [V_{h}]^2,\\
	&\int_{\Omega}\nabla\cdot \bm{w}_k^nq_Hdx=0, \forall q_H\in P_H.
	\end{aligned}
	\right.
	\end{flalign*}
Consequently, the following quadratic function
\begin{equation*}%\label{objective_function_step size}
	Q_k(\rho)=\frac{1}{2}{\Delta t}\sum_{n=1}^N\int_{\Omega} |\rho\bm{w}_k^n +\bm{y}_k^n-y_d^n|^2dx+\frac{\alpha}{2}\Delta t \sum_{n=1}^N \int_{\Gamma_c} |\bm{u}^n_k+\rho \bm{d}^n_k|^2dx
\end{equation*}
is an approximation of $P_k(\rho)$. Then it is easy to show that the minimizer of $Q_k(\rho)$ is
\begin{equation}\label{step_bfgs}
\hat{\rho}_k=-\frac{{\Delta t}\sum_{n=1}^N\int_{\Omega}  \bm{g}_k^n \cdot \bm{d}_k^n dx}{{\Delta t}\sum_{n=1}^N\int_{\Omega} |\bm{w}_k^n|^2dx+\alpha \Delta t \sum_{n=1}^N \int_{\Gamma_c} | \bm{d}^n_k|^2dx
},
\end{equation}
where $(\bm{g}_k^n)_{n=1}^N=DJ_h^{\Delta t}(\bm{u}_k)$.
We take step size $\rho_k=\hat{\rho_k}$, which is indeed an approximate minimizer of $P_k(\rho)$.

\subsection{Computation of the search direction ${d}_k$}\label{se:bfgs_gradient}
As mentioned in Section \ref{se:generic_bfgs}, due to the requirement of large memory, it is more appropriate to consider the L-BFGS method, which only requires storing a sequence of vectors $\delta \bm{u}_k$ and $\delta \bm{g}_k$ from the most recent iterations to compute $\mathcal{H}_k\bm{g}_k$ without constructing $\mathcal{H}_k$ explicitly. It turns out that the L-BFGS method is fairly robust, inexpensive, and easy-to-implement, see e.g., \cite{nocedal1980,nocedal2006sequential}.
For the reader's convenience, we present a two-loop recursive procedure to efficiently compute the product $\mathcal{H}_k\bm{g}_k$ in Algorithm \ref{alg:two_loop}.

\begin{algorithm}
	\caption{Computation of search direction}
	\label{alg:two_loop}
	\begin{algorithmic}
		\STATE Set $\hat{\bm{d}}=-\bm{g}_k (= -DJ_h^{\Delta t}(\bm{u}^k))$;
		\FOR {$i=k-1,\cdots,k-m$}
		\STATE$\eta_i=\frac{1}{\delta \bm{g}_i^T \delta \bm{u}_i}$, $\tau_i=\eta_i\delta \bm{u}_i^T \hat{\bm{d}}$
		\STATE$\hat{\bm{d}}=\hat{\bm{d}}-\tau_i\delta \bm{g}_i$
		\ENDFOR
		\STATE Compute $\bm{d}=\mathcal{H}_k^0\hat{\bm{d}}$
		\FOR {$i=k-m, k-m+1, \cdots, k-1$}
		\STATE $\kappa=\eta_i \delta \bm{g}_i^T \bm{d}$
		\STATE $\bm{d}=\bm{d}+\delta \bm{u}_i(\tau_i-\kappa)$
		\ENDFOR
		\RETURN $\bm{d}=-\mathcal{H}_k\bm{g}_k$
	\end{algorithmic}
\end{algorithm}

Note that, in contrast to the standard BFGS iteration,
the initial inverse Hessian approximation $\mathcal{H}_k^0$ is allowed to vary from iteration to iteration.
In our numerical implementation, we choose $\mathcal{H}_k^0=I$ for simplicity.
% For instance, one can choose
%$\mathcal{H}_k^0=\gamma_k\mathcal{H}_0^0, \quad \text{with}\quad \gamma_k>0~ \text{and}~ \mathcal{H}_0^0=I.$

\subsection{An easily implementable L-BFGS method for the solution of problem (\ref{fully_discrete_optimal_control})-(\ref{fully-discretized-state-variationalform})}
With the discussions in Sections \ref{se:bfgs_step} and \ref{se:bfgs_gradient}, we obtain an easily implementable L-BFGS method as listed in Algorithm \ref{alg:L_CBFGS} for the solution of problem (\ref{fully_discrete_optimal_control})-(\ref{fully-discretized-state-variationalform}).
\begin{algorithm}
	\caption{An easily implementable L-BFGS method for the solution of problem (\ref{fully_discrete_optimal_control})-(\ref{fully-discretized-state-variationalform})}
	\label{alg:L_CBFGS}
	\begin{algorithmic}
		\STATE Given an initial point $\bm{u}_0$.
		\STATE Compute $\bm{g}_0=D J_h^{\Delta t}(\bm{u}_0)$; If $\bm{g}_0=0$, then $\bm{u}=\bm{u}_0$, else for $k\geq 0$
	\WHILE {not converged}
		\STATE Choose $\mathcal{H}_k^0$;
		\STATE Compute search direction $\bm{d}_k$ by Algorithm \ref{alg:two_loop}:
		\STATE Compute a step size $\rho_k$ by (\ref{step_bfgs}).
		\STATE Update $\bm{u}_{k+1}=\bm{u}_k+\rho_k \bm{d}_k$ and $\bm{g}_{k+1}=D J_h^{\Delta t}(\bm{u}_{k+1})$
		\IF {$k>m$}
		\STATE Discard the vector pair $\{\delta \bm{u}_{k-m}, \delta \bm{g}_{k-m}\}$ from storage
		\ENDIF
		\STATE Compute and save $\delta \bm{u}_k=\bm{u}_{k+1}-\bm{u}_k$ and $\delta \bm{g}_k=\bm{g}_{k+1}-\bm{g}_k$
	\ENDWHILE
	\end{algorithmic}
\end{algorithm}

\section{Numerical results}\label{sec:numerical_results}

In this section, we report some preliminary numerical results to validate the efficiency of our proposed numerical approach for solving the optimal control problem (\ref{OptimalControl}) with the objective function (\ref{objective_functional}) or (\ref{objective_functional_r}).  In particular, we test two numerical examples for velocity-tracking and vorticity reduction by controlling the heat flux. All codes were written in MATLAB R2020b and numerical experiments were conducted on a
MacBook Pro with mac OS Monterey, Intel(R) Core(TM) i7-9570h (2.60 GHz), and 16 GB RAM.
	 For simplicity, all the linear equations originated from the numerical discretization are solved by the default backslash operator in Matlab.

%	 	The difference is we use $\bm{y}_h^n \in V_{0h}$
%	 	via solving
%	 	\begin{flalign}\label{deg_Stokes_y_2_dis}
%	 		(\frac{\bm{y}_h^{n}-\tilde{\bm{y}}_h^n}{\Delta t}, {\varphi_h})-b(\bm{p}_H^n-\bm{p}_H^{n-1},{\varphi_h})=0,\quad \forall {\varphi_h}\in [V_{0h}]^2.
%	 	\end{flalign}
%	 	which is the fully discrete analogue of

	  % After space-time discretization, a finite-dimensional nonlinear optimization problem will need to be solved.
	
	\medskip
\noindent\textbf{Example 1.}  We consider the tracking-type optimal control of the Boussinesq equations (\ref{state_equation1})-(\ref{state_equation2}) with the objective functional \eqref{objective_functional}. The boundary and initial conditions are specified as
\begin{flalign*}%\label{state_equation_E1}
\begin{aligned}
&y=0~ \text{on} ~\Gamma\times(0,T),\quad
&y(0)=0.&\\
&\left\{\begin{aligned} &\theta=0&~ \text{on} ~\Gamma\slash\Gamma_c \times(0,T),\\ &\nu_2\frac{\partial\theta}{\partial \vec{n}}=v &~\text{on} ~\Gamma_c \times(0,T),\end{aligned}\right. \quad
&\theta(0)=0.
\end{aligned}
\end{flalign*}

\begin{figure}[htpb]
	\centering{\includegraphics[width=0.45\textwidth]{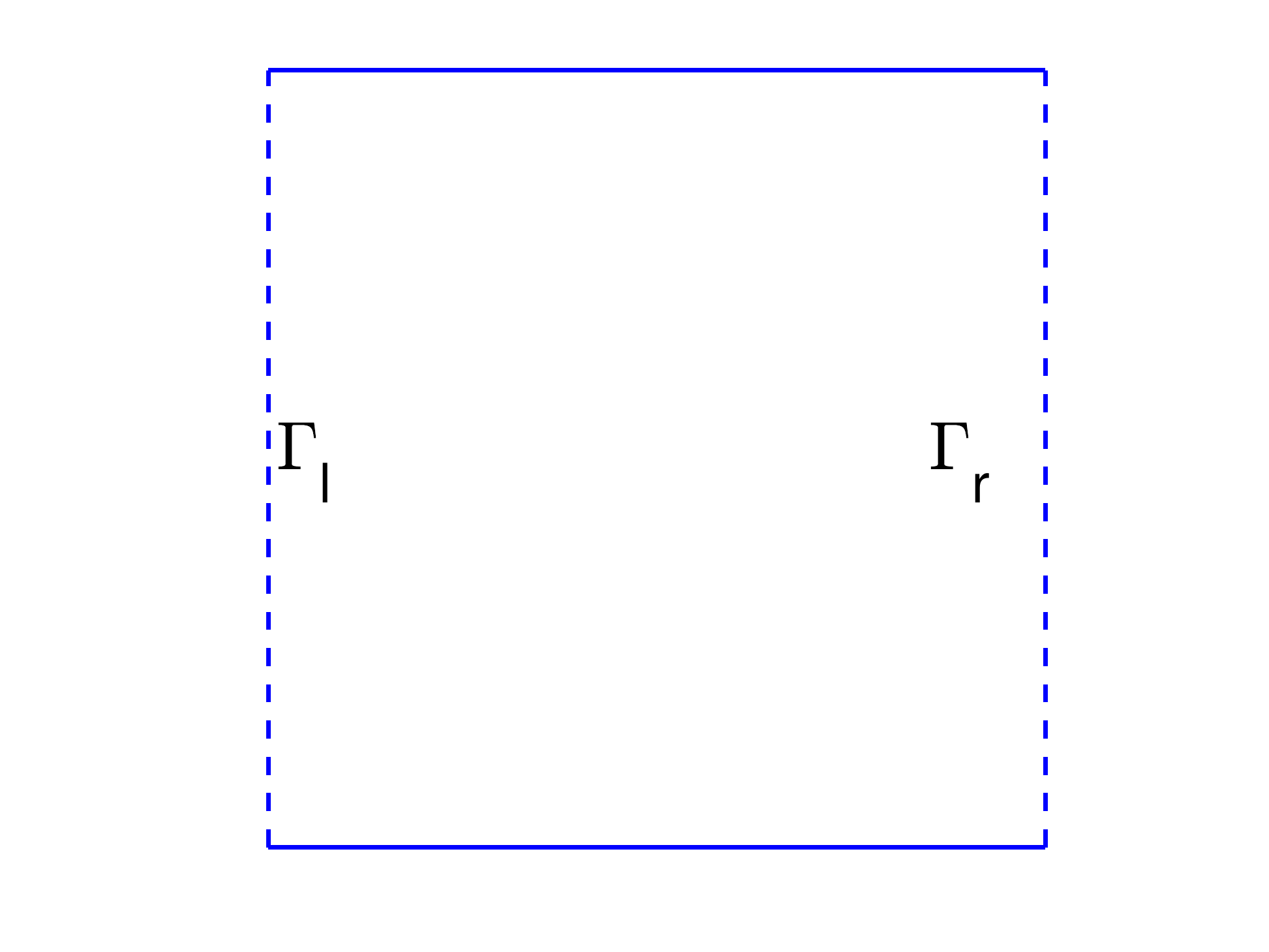}}
	\centering{\includegraphics[width=0.45\textwidth]{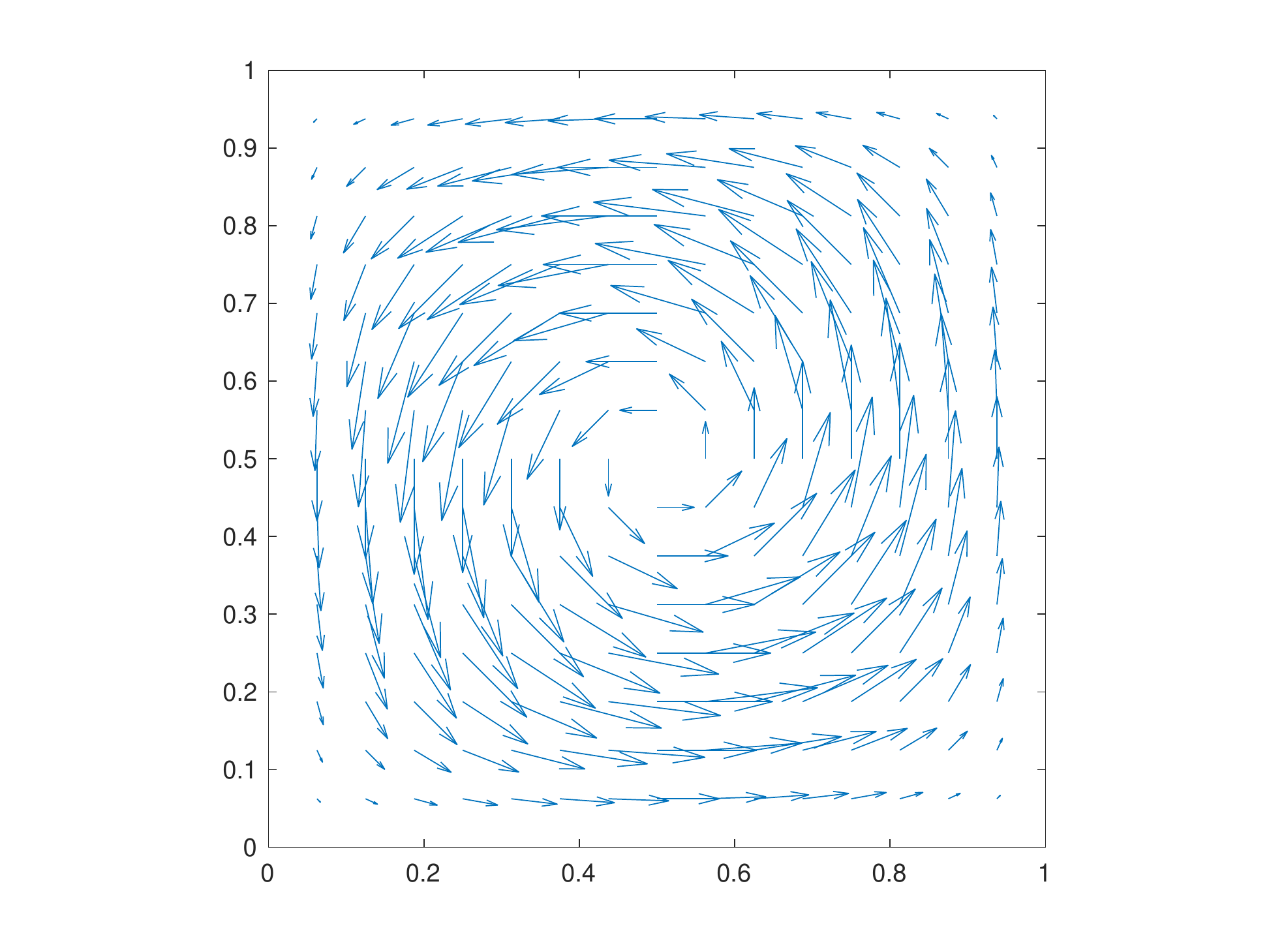}}
	\caption{The schematics of $\Gamma_c=\Gamma_{l}\cup\Gamma_{r}$ (left) and the target velocity ${y}_d$ (right) for Example 1..}
	\label{Domain_1}
\end{figure}

Furthermore, we set that
\begin{itemize}
	\item the regularization parameter is $\alpha=5\times10^{-5}$ and the coefficients are $\nu_1=1/100$ and $\nu_2=1/72$;
	\item the domain is $\Omega=(0,1)^2$, $T=5$, and the control acts on $\Gamma_c=\Gamma_{l}\cup\Gamma_{r}$  as shown in Figure \ref{Domain_1} (left);
	\item the target velocity ${y}_d=(y_{1,d},y_{2,d})^T$ is time-independent and given by
$$
\left\{
\begin{aligned}
y_{1,d}(x_1,x_2,t)&=100\phi(x_1)\phi'(x_2),\\
y_{2,d}(x_1,x_2,t)&=-100\phi'(x_1)\phi(x_2),\\
\end{aligned}
\quad\forall~ (x_1,x_2)\in \Omega, t\in (0,T)
\right.
$$
where $\phi(z)=z^2(z-1)^2$ and $\nabla \cdot {y}_d=0$.  A schematic of the target velocity ${y}_d$ is shown in Figure \ref{Domain_1} (right).
\end{itemize}

If there is no control in Example 1, i.e., the control is set as $v=0$, then the states $y$ and $\theta$ are zeros in $(0, T)$.  We aim at tracking a time-independent pinwheel-shape velocity field $y_d$ from static initial states $y(0)=0$ and $\theta(0)=0$  by controlling the heat flux $v$ on the side walls $\Gamma_c$ of the square domain $\Omega$. For this purpose, one can generate a buoyancy-driven swirling flow whose velocity field $y$ is close to the target field $y_d$ during $(0, T)$.

%To this end, we heat the upper half of the left side wall and cool the lower half of the right side wall, then we will get a swirling flow driven by buoyancy, and its velocity field $y$ is desired to be close to the target field $y_d$.

\begin{figure}[htpb]
	\centering{\includegraphics[width=0.45\textwidth]{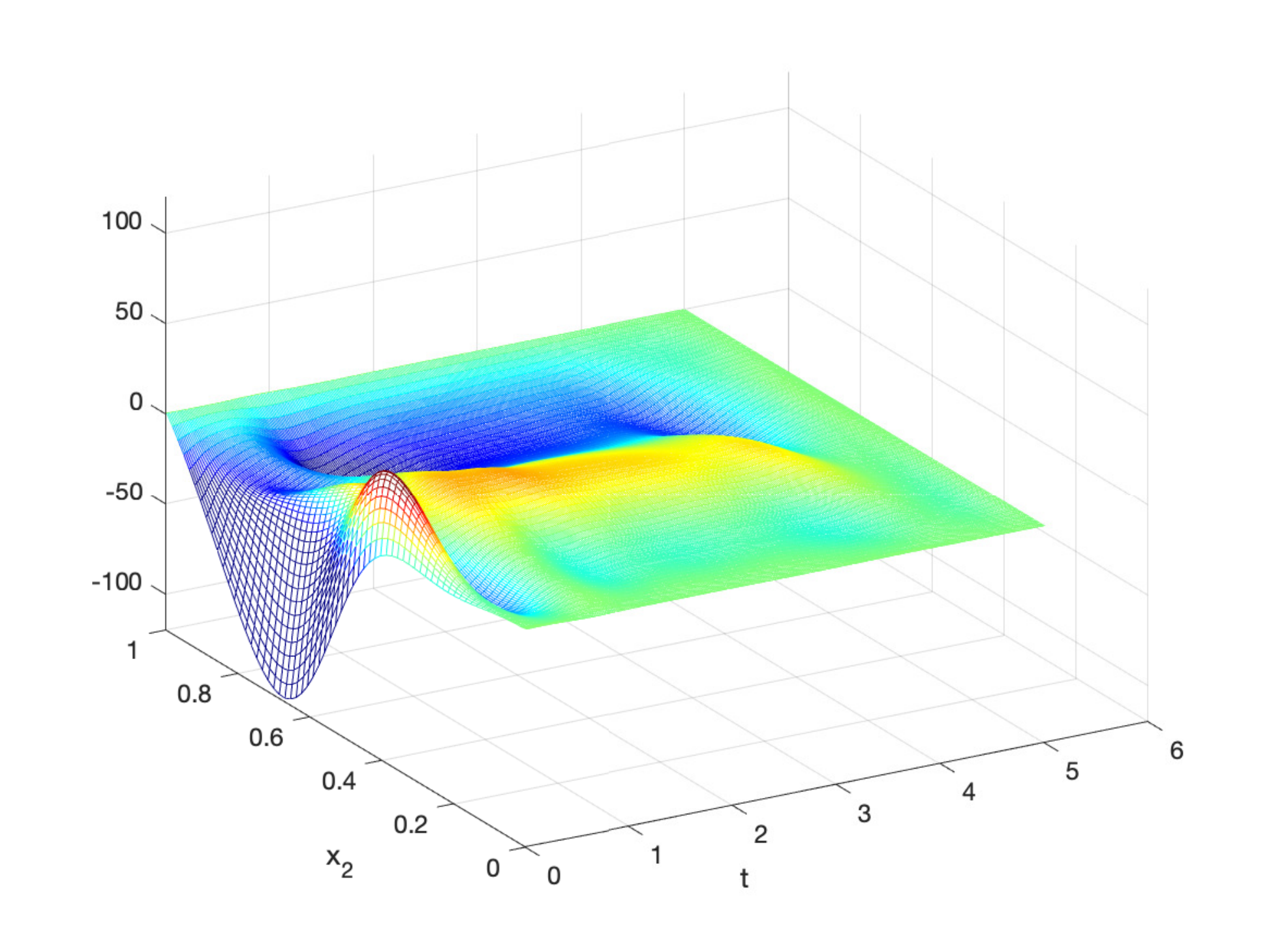}}
	{\includegraphics[width=0.45\textwidth]{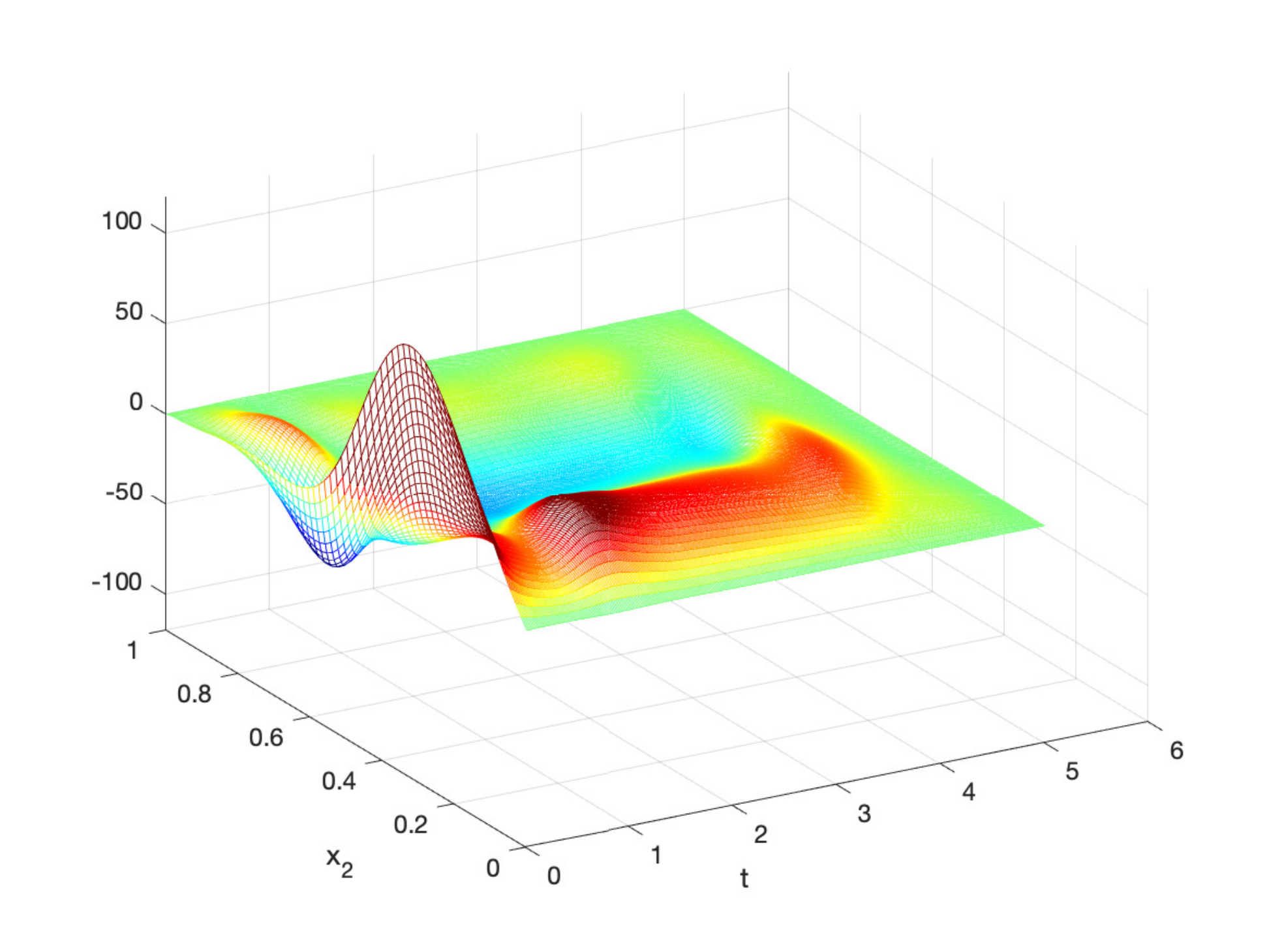}}
	\caption{The computed optimal control $u$ on  $\Gamma_{l}$ (left) and  $\Gamma_{r}$ (right) for Example 1.}
	\label{Control_1}
\end{figure}

\begin{figure}[htpb]
	\centering{\includegraphics[width=0.45\textwidth]{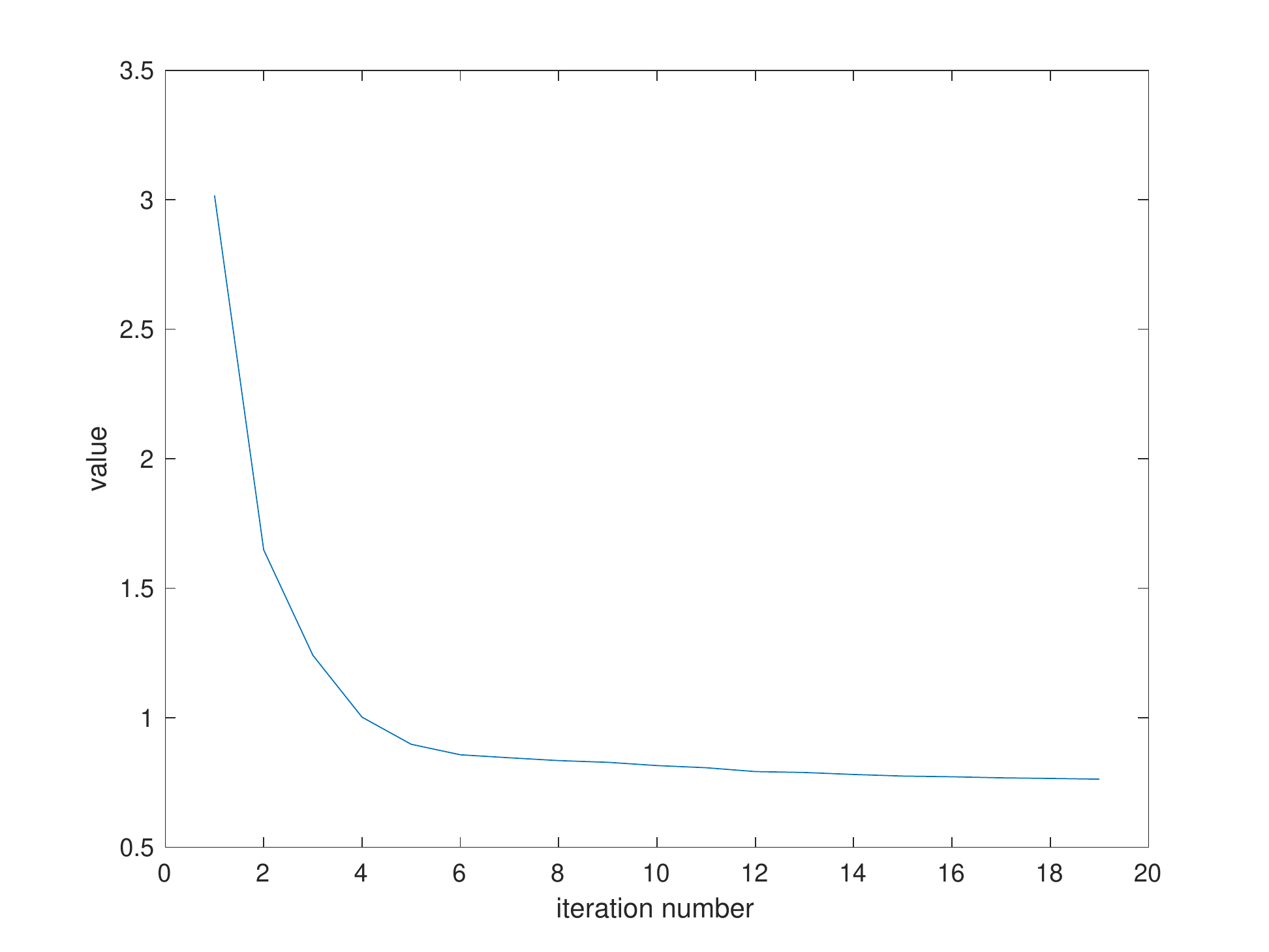}}
	{\includegraphics[width=0.45\textwidth]{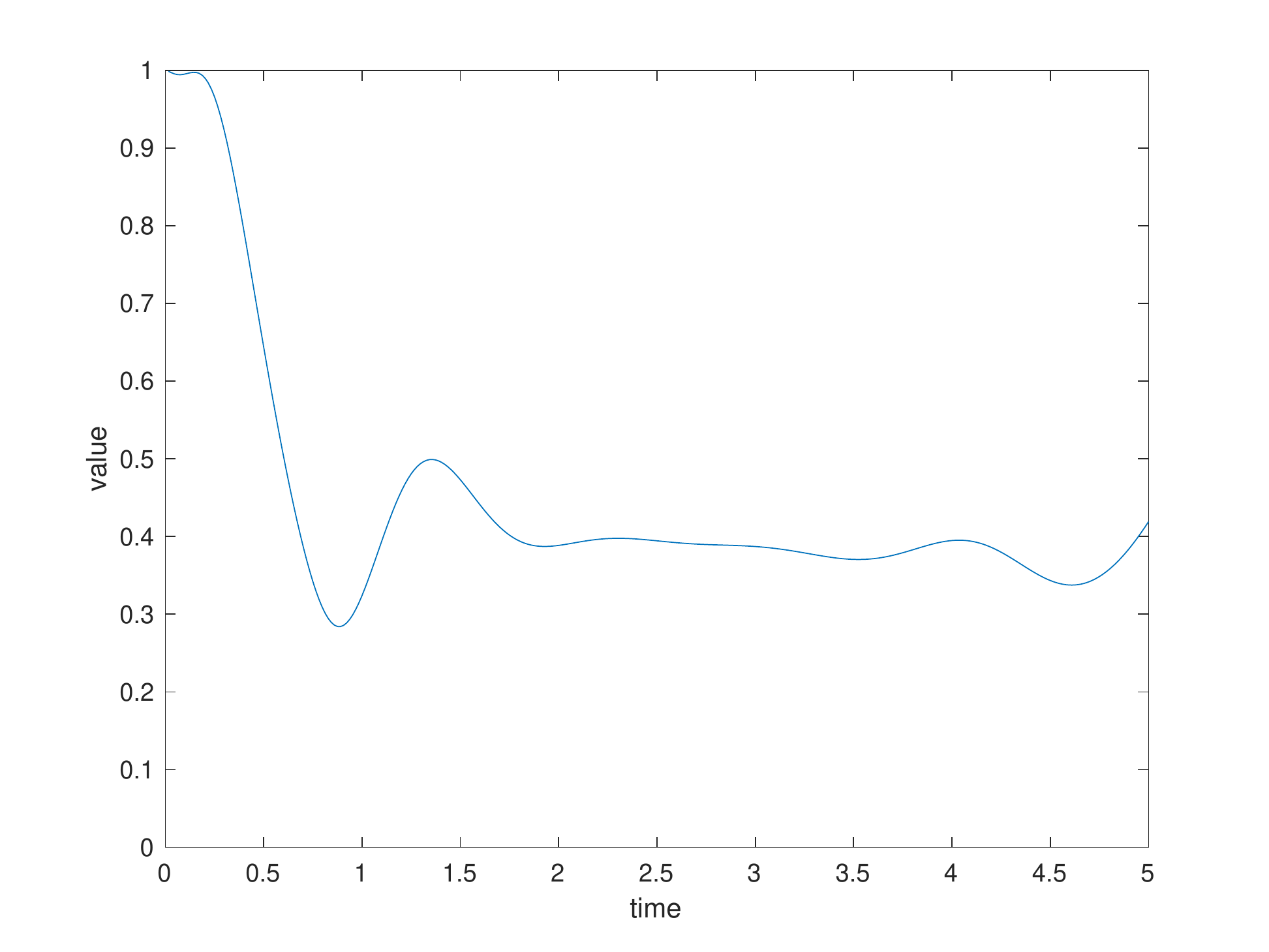}}
	\caption{Numerical results for Example 1. Left: the objective functional values with respect to L-BFGS iterations.
	Right: relative errors $\frac{\|y(t;u)-y_d(t)\|}{\|y_d(t)\|}$ at different time $t$.
		}
	\label{obj_1}
\end{figure}

\begin{figure}[htpb]
	\centering
	{\includegraphics[width=0.32\textwidth]{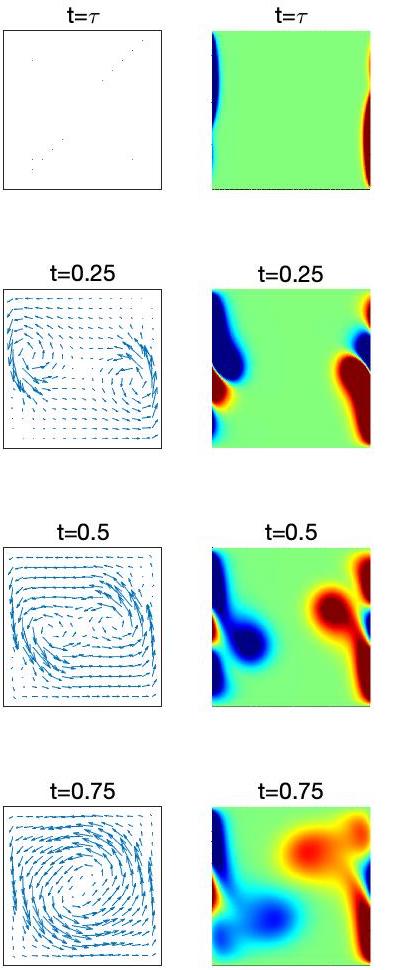}}
	{\includegraphics[width=0.32\textwidth]{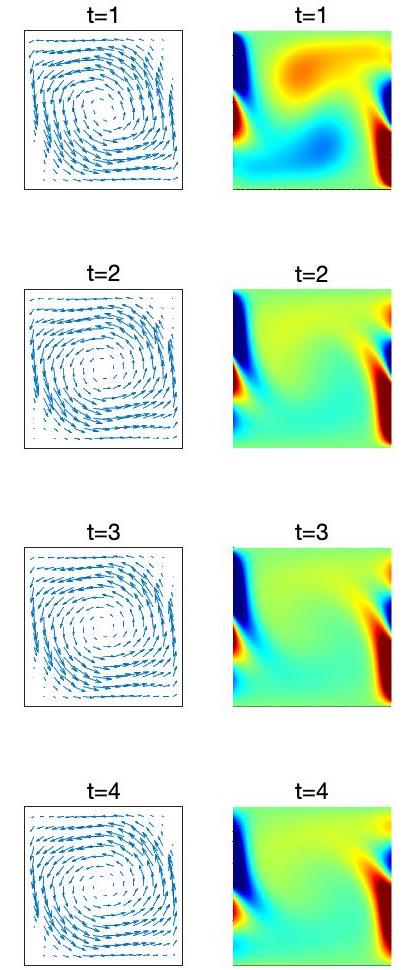}}
	{\includegraphics[width=0.32\textwidth]{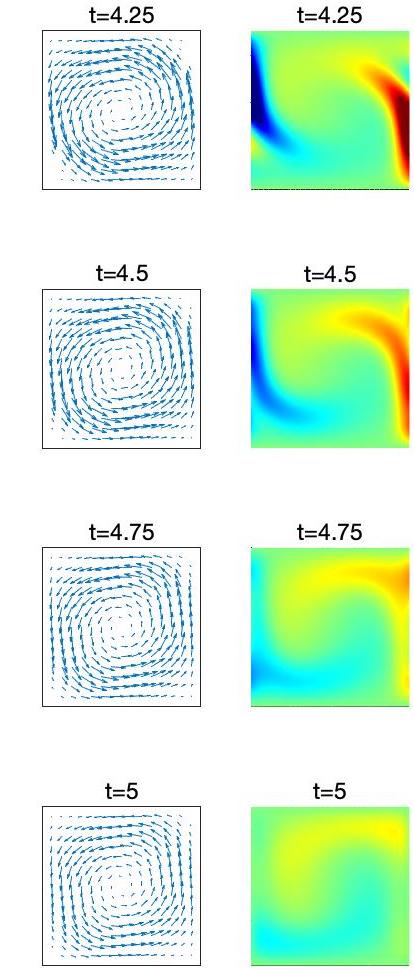}}
	\caption{Computed velocity field $y$ (the 1st, 3rd, and 5th columns) and  computed temperature $\theta$ ((the 2nd, 4th, and 6th columns) ) at different time $t$ for Example 1. Here,  blue/red color means low/high temperature. }
	\label{States_1}
\end{figure}

We implement Algorithm \ref{alg:L_CBFGS} to Example 1 with mesh size $h=1/64$ for the space discretization and step size $\Delta t=1/64$ for the time discretization.  Moreover, we store the most recent vector pair $\{\delta \bm{u}_{k-m}, \delta \bm{g}_{k-m}\}$ during the L-BFGS iterations and terminate it when the following stopping criterion is satisfied:
$$
\frac{\|DJ(\bm{u}_h^k)\|}{\|D J(\bm{u}_h^0)\|}<5\times10^{-3}.
$$
The numerical results for Example 1 are displayed in Figures \ref{Control_1},  \ref{obj_1}, and \ref{States_1}.

%The computed optimal controls on the left and right sides are anti-symmetric at any time, which is reasonable because the target field $y_d$ is anti-centrosymmetric.
We observe from Figure \ref{obj_1} that the objective function values decrease rapidly, which implies fast convergence and hence high efficiency of our proposed Algorithm \ref{alg:L_CBFGS}.  The relative error $\|y(t;u)-y_d(t)\|/\|y_d(t)\|$ becomes small after about one second, which indicates that the velocity field gets close to the target field, and this observation coincides with that in Figure \ref{States_1}.

% {\color{blue} In addition, we note that the norm of the computed optimal control $\|u_h(t)\|$ is large initially, then decreases rapidly, and finally goes to zero. From these observations, we can conclude that it costs lots of energy to drive the initial state $y(0)$ to be close to the target state $y_d$ at the beginning of the control process. Then, the cost decreases when the state is kept unchanged, and finally, the control is gradually removed.}

The computed velocity and temperature at different instants of time for Example 1 are displayed in Figure \ref{States_1}, from which we can observe how the buoyancy drives a swirling fluid and how the velocity influences the temperature.  To be concrete, at the early stage of the control process ($t=0.016,0.25,0.5,0.75$), the temperature and the velocity field begin to change under the action of the optimal control. Two swirling flows are formed on the left and right sides and merge into a large swirling flow.
 Then, at $t=1, 2, 3, 4$, both the velocity and the temperature become steady and the velocity field is close to the target field $y_d$.  In the final stage $(t=4.25, 4.5, 4.75, 5)$, the control is gradually removed and the temperature gradually goes to 0, but the computed velocity field remains close to the target field $y_d$.

\medskip

\noindent\textbf{Example 2.}  This example is extended from \cite{IR98}, which is motivated by the transport process in high-pressure chemical vapor deposition (CVD) reactors. The original example focuses on the steady case and we extend it to an unsteady example in $(0, T)$ with $T=15$.

\begin{figure}[htpb]
	\centering{\includegraphics[width=0.45\textwidth]{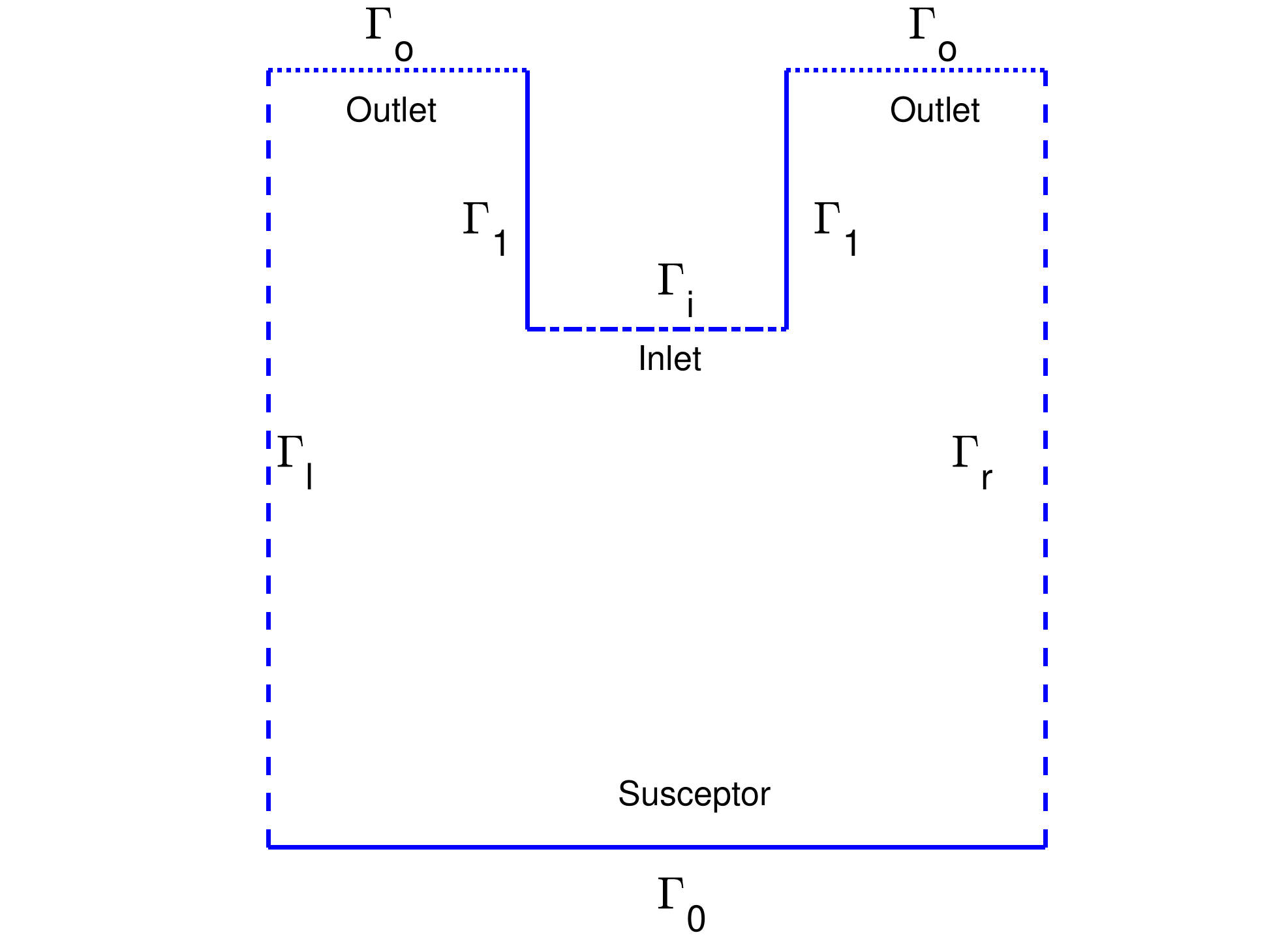}}
	\centering{\includegraphics[width=0.45\textwidth]{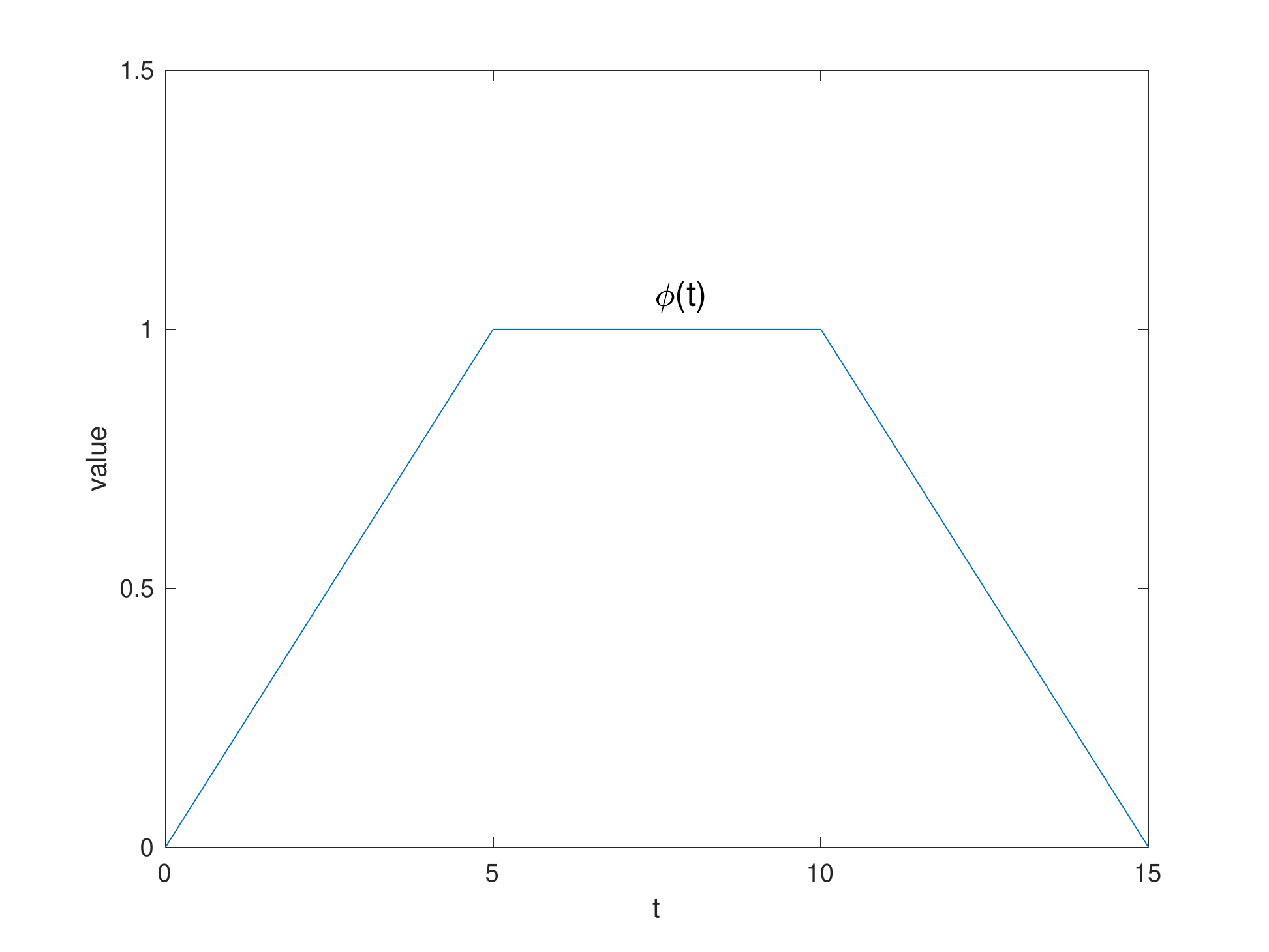}}
	\caption{The schematics of the domain $\Omega$ (left) and the graph of $\phi(t)$ (right) for Example 2. a) The reactant gases entered the reactor from $\Gamma_i$, then flow down to the substrate $\Gamma_0$ which is kept at high temperature, and finally go out from $\Gamma_o$. b) The function $\phi(t)$ determines the maximum velocity of the fluid on $\Gamma_i$. It increases linearly from zero to one during $(0,5)$, then remains unchanged during $(5, 10)$, and finally decreases linearly to zero during $(10, 15)$. }
	\label{Domain_2}
\end{figure}

 A typical vertical reactor $\Omega$, shown in Figure  \ref{Domain_2}, is a classical configuration for the growth of compound semiconductors by metalorganic vapor phase epitaxy. The geometry of the prototype reactors has two outlet portions, $\Gamma_o$, and an inlet, $\Gamma_i$, whose widths are 1/3.
The size of the susceptor region $\Gamma_0$ and that of the side walls $\Gamma_{l}$ and $\Gamma_{r}$ are 1; the height of the inlet port $\Gamma_1$ is 1/3.  The reactant gases enter the reactor from $\Gamma_i$ and flow down to the substrate $\Gamma_0$ which is kept at a high temperature. This means that the least dense gas is closest to the substrate and the flow is likely to be affected by buoyancy-driven convection. To have uniform growth rates and better compositional variations, it is crucial to have a flow field without recirculation.

To achieve the aforementioned goal, one can minimize the vorticity by controlling the temperature at the side walls $\Gamma_c=\Gamma_{l}\cup\Gamma_{r}$ in order to obtain a flow field without recirculation and thereby obtain better vertical transport. For this purpose, we consider
to minimize the  objective functional
$$
J(v)=\frac{1}{2}\iint_Q|\nabla \times y |^2dxdt+\frac{\alpha}{2}\iint_{\Sigma_c} v^2dxdt,
$$
where $y:=y(v)$ is the solution of the state equation \eqref{state_equation1} complemented with the following initial and boundary conditions:
\begin{flalign}\label{state_equation_E2}
\left\{
\begin{aligned}
&y(0)=0, \theta(0)=0\\
&\Gamma_0: y=(0,0),\quad \theta=1,\\
&\Gamma_1: y=(0,0), \quad \theta=0,\\
&\Gamma_i: y=(0,-4(x_1-1/3)(2/3-x_1)\phi(t)), \quad \theta=0,\\
&\Gamma_0: \frac{\partial y}{\partial \vec{n}}=(0,0), \quad \frac{\partial \theta}{\partial \vec{n}}=0,\\
&\Gamma_c, y=(0,0), \quad \frac{\partial \theta}{\partial \vec{n}}+\theta=v.
\end{aligned}
\right.
\end{flalign}
Here the function
$$\phi(t):=
\left\{\begin{aligned}
	&t/5,\quad& 0<t\le 5,\\
	&1,\quad& 5<t<10,\\
	&(15-t)/5,\quad& 10\le t<15,
\end{aligned}
\right.
$$ determines the maximum velocity of the fluid on $\Gamma_i$.

Throughout, the regularization parameter is $\alpha=10^{-4}$ and the coefficients are $\nu_1=1/100$ and $\nu_2=1/72$ respectively.
We use a uniform triangulation with mesh size $h=1/64$ for the space discretization and step size $\Delta t=1/64$ for the time discretization.
We keep the $m=5$ most recent vector pair $\{\delta \bm{u}_{k-i}, \delta \bm{g}_{k-i}\}_{i=1}^m$ during the L-BFGS iterations and terminate it if
$
\frac{\|D J(\bm{u}_h^k)\|}{\|DJ(\bm{u}_h^0)\|}<5\times10^{-3}.
$
\begin{figure}[htpb]
	\centering{\includegraphics[width=0.45\textwidth]{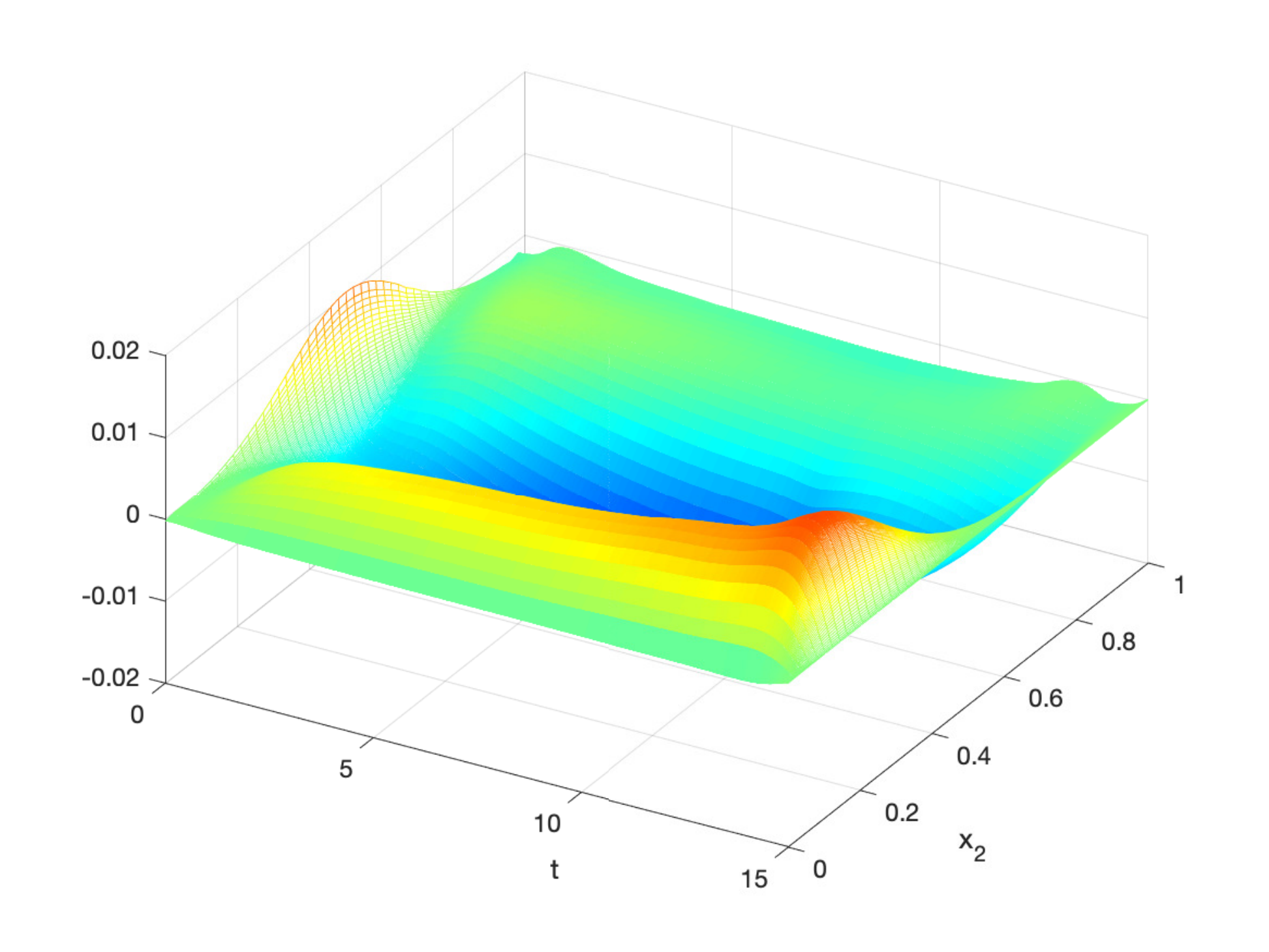}}
	{\includegraphics[width=0.45\textwidth]{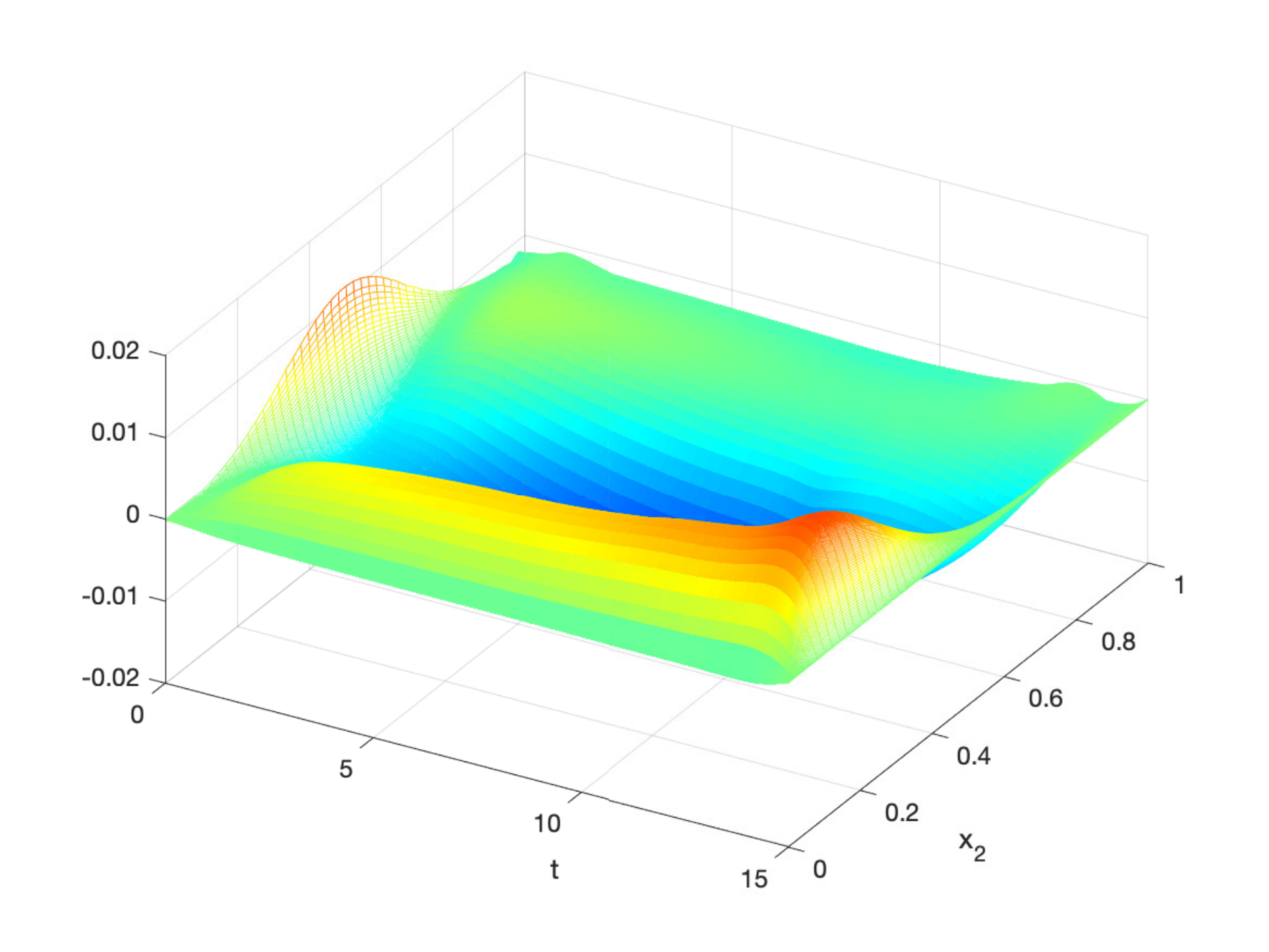}}
	\caption{The computed optimal control $u$ on  $\Gamma_{l}$ (left) and  $\Gamma_{r}$ (right) for Example 2.
		}
	\label{Control_2}
\end{figure}
\begin{figure}[htpb]
	\centering{\includegraphics[width=0.45\textwidth]{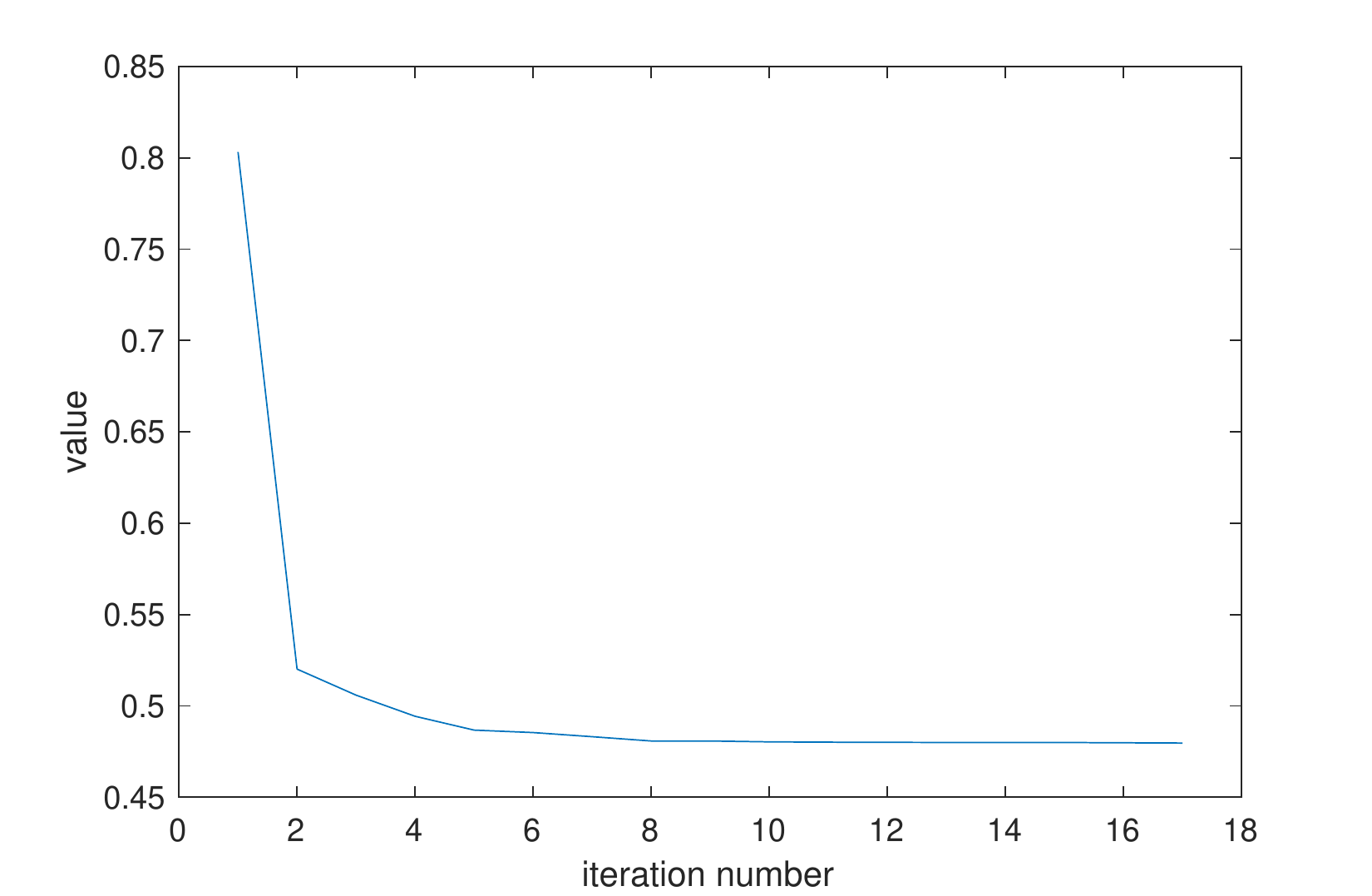}}
	{\includegraphics[width=0.45\textwidth]{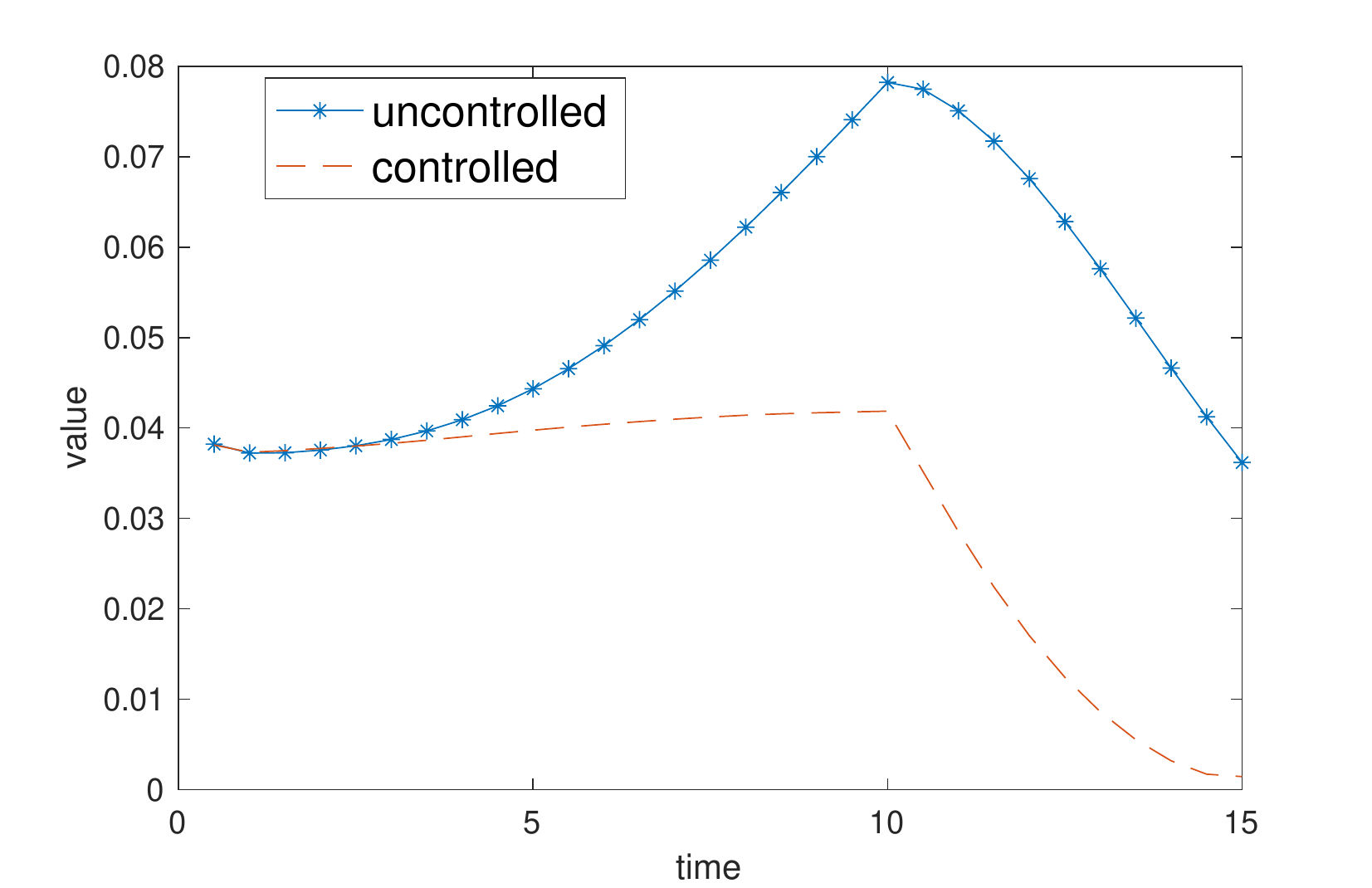}}
	\caption{Numerical results for Example 2. Left:  the objective functional values with respect to L-BFGS iterations. Right:  the vorticity values $\|\nabla \times y(t;u)\|$ at different times t.
		}
	\label{obj_2}
\end{figure}

\begin{figure}[htpb]
	\centering{\includegraphics[width=0.45\textwidth]{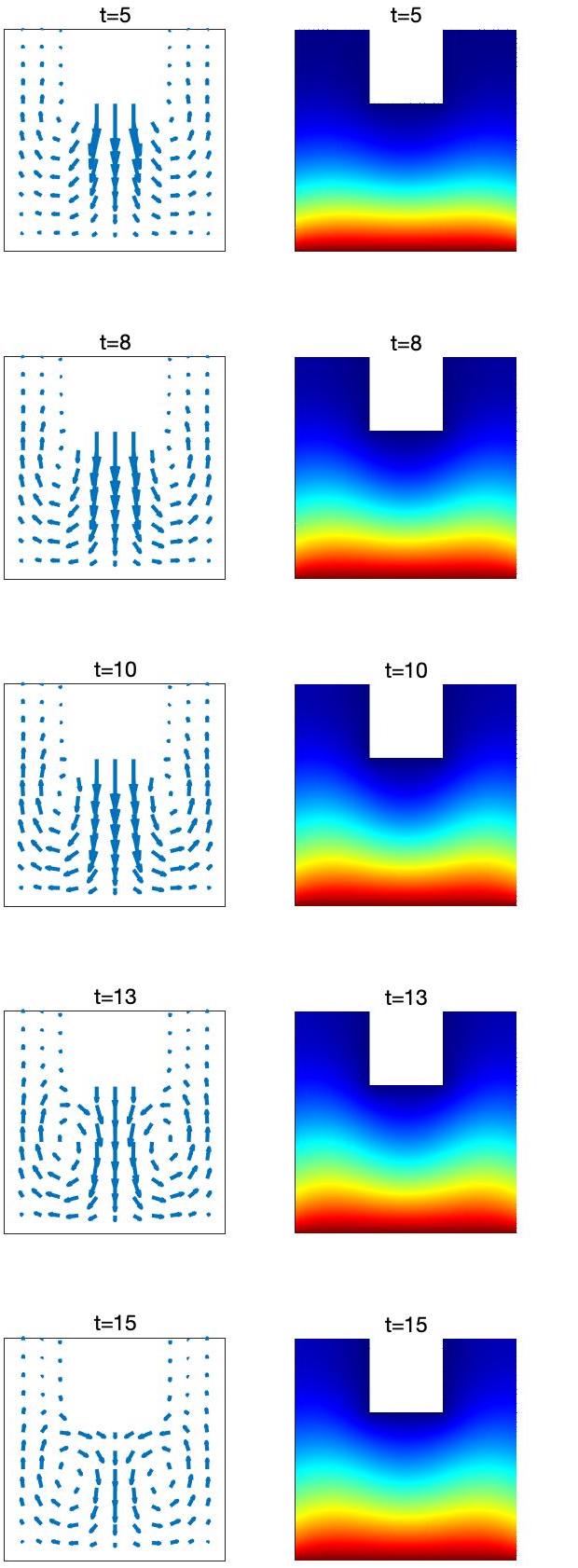}}
	\centering{\includegraphics[width=0.45\textwidth]{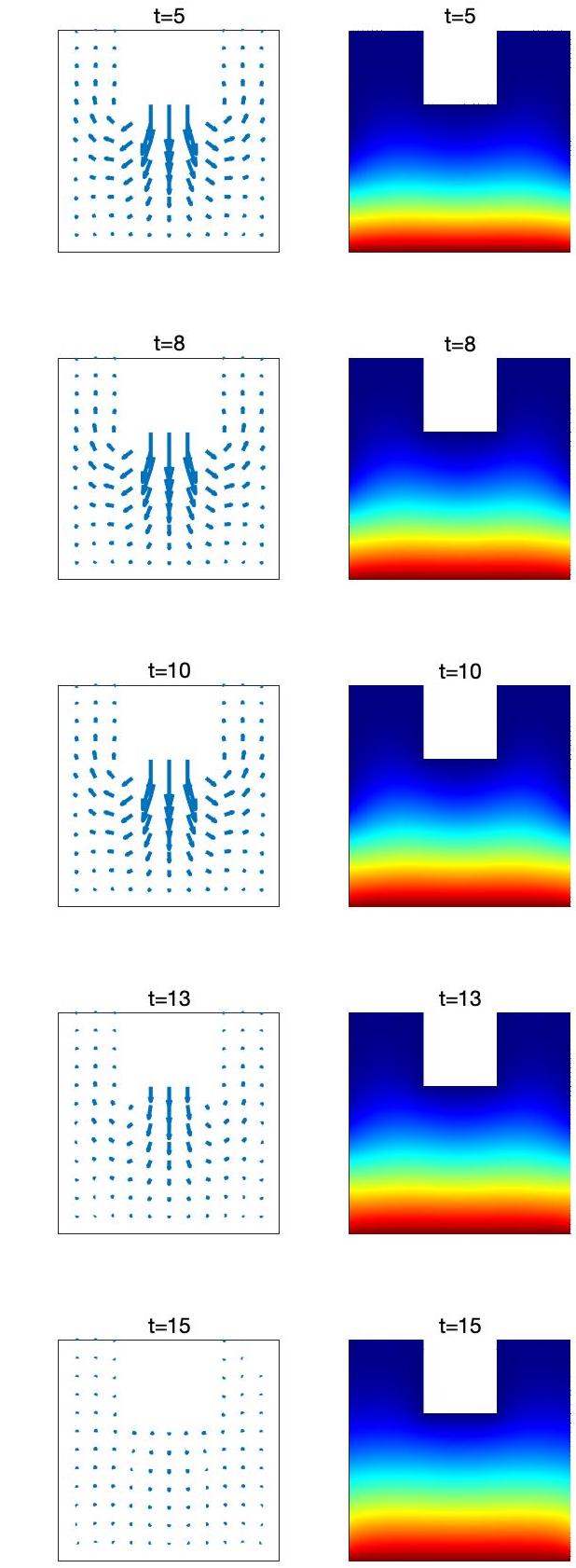}}
	\caption{Computed velocity field $y$ (the 1st and 3rd columns) and temperature $\theta$ (the 2nd and 4th columns) of the uncontrolled (left) and controlled (right) system at different $t$ for Example 2. }
\label{States_2}
\end{figure}
The numerical results obtained by Algorithm \ref{alg:L_CBFGS} for Example 2 are displayed in Figures \ref{Control_2}, \ref{obj_2}, and \ref{States_2}, where the notation "uncontrolled" means that there is no control (i.e. $v=0$) acting in the state system specified by \eqref{state_equation1} and \eqref{state_equation_E2}.
  We observe from Figure \ref{Control_2} that the computed controls on the left and right sides are the same at any time, which is due to the symmetry structure of the example under investigation.   From Figure \ref{obj_2}, we see that the objective function values decrease rapidly which implies fast convergence of Algorithm \ref{alg:L_CBFGS}. Moreover, we observe that the vorticity  $\|\nabla\times y(t,u)\|$ of the controlled system goes to zeros eventually. We observe from Figure \ref{States_2}  that swirling flow appears near the susceptor for the uncontrolled system while there is no swirling flow for the controlled system; and the difference is more discernible after $t=10$. This implies that the computed control works very well and it indeed leads to a flow field without recirculation. This, together with the results shown in Figure \ref{obj_2}, validates a significant vorticity reduction by the computed control.

 \section{Conclusions}
\label{sec:conclusions}

In this paper, we proposed an efficient numerical approach to the optimal control of thermally convective flows, which can be mathematically modeled as optimally controlling the Boussinesq equations consisting of the Navier-Stokes equations for incompressible viscous flow coupled with an advection-diffusion equation for temperature. Our main techniques included the Marchuk-Yanenko operator splitting method for the time discretization to untie the advection-diffusion equation from the Navier-Stokes component and to decouple the nonlinearity from the incompressibility condition. Computing the gradient of the objective functional became possible, and it was reduced to solving four easy linear advection-diffusion equations and two degenerated Stokes equations at each time step. With the Bercovier-Pironneau finite element method for space discretization, we also proposed an efficient and easily implementable BFGS method with limited memory to solve the fully discretized optimal control problem. The proposed algorithm appears to be the first efficient numerical approach to the optimal control of unsteady thermally convective flows.

We focused on the numerical study for optimal control problems modeled by the Boussinesq equations, and its validated efficiency clearly justifies the necessity to investigate the underlying theoretical issues such as the convergence analysis for the proposed numerical approach and the error estimate for the operator splitting time discretization scheme. Moreover, we conjecture that our philosophies in algorithmic design and techniques for numerical implementation can be extended to other important optimal control problems modeled by coupled systems, such as the optimal mixing \cite{HZ2021}, the optimal control of coupled  Cahn-Hilliard Navier-Stokes system \cite{HW2014}, and the optimal control of a diffuse interface model for tumor growth \cite{ebenbeck2020}.

\section*{Acknowledgments} The authors would like to thank Prof. Ming-Chih Lai for his valuable suggestions about the numerical experiments.

\end{document}